\pdfoutput=1
\documentclass[12pt]{amsart}
\usepackage{amssymb}
\usepackage{cite}
\usepackage[margin=1in,dvips]{geometry}
\usepackage{booktabs}
\usepackage{url}
\usepackage{hyphenat}
\usepackage[bookmarks=true,%
    colorlinks=true,%
    linkcolor=blue,%
    citecolor=blue,%
    filecolor=blue,%
    menucolor=blue,%
    pagecolor=blue,%
    urlcolor=blue]{hyperref}
\usepackage{mathtools}
\usepackage{tikz}
\usetikzlibrary{matrix,arrows}
\tikzset{cdarrow/.style={auto,
    execute at begin node=$\scriptstyle,execute at end node=$}}

\newread\testin

\def\mathcenter#1{%
  \vcenter{\hbox{$#1$}}%
}

\def\graph#1{
        \includegraphics{#1}
}

\def\mathgraph#1{
        \mathcenter{\graph{#1}}
}

\graphicspath{{draws/}{mpdraws/}{picts/}}

\newcommand{\RR}{\mathbb R}
\newcommand{\CC}{\mathbb C}
\newcommand{\ZZ}{\mathbb Z}
\newcommand{\NN}{\mathbb N}
\newcommand{\QQ}{\mathbb Q}
\newcommand{\PP}{\mathbb P}
\newcommand{\EE}{\mathbb E}

\newcommand{\abs}[1]{{\lvert #1 \rvert}}

\newcommand{\isom}{\simeq}

\newcommand{\superset}{\supset}

\newcommand{\tensor}{\otimes}

\DeclareMathOperator{\rank}{rank}
\DeclareMathOperator{\af}{Aff}
\DeclareMathOperator{\ls}{\spanop}
\DeclareMathOperator{\ann}{ann}
\DeclareMathOperator{\Gr}{Gr}

\DeclareMathOperator{\spanop}{span}
\DeclareMathOperator{\chord}{chord} % 
\DeclareMathOperator{\adj}{adj} % adjoint matrix

\DeclareMathOperator{\Eucl}{Eucl}
\DeclareMathOperator{\Trans}{Trans}

\theoremstyle{plain}
\newtheorem{proposition}{Proposition}[section]
\newtheorem{theorem}[proposition]{Theorem}
\newtheorem{lemma}[proposition]{Lemma}
\newtheorem{corollary}[proposition]{Corollary}

\theoremstyle{definition}
\newtheorem{definition}[proposition]{Definition}

\newtheorem{algorithm}[proposition]{Algorithm}

\theoremstyle{remark}

\newtheorem{remark}[proposition]{Remark}
\newcommand{\Edges}{{\mathcal E}}
\newcommand{\Verts}{{\mathcal V}}

\newcommand{\tA}{\widetilde A}
\newcommand{\tB}{\widetilde B}
\newcommand{\tl}{\widetilde \ell}
\newcommand{\trh}{\tilde{\rho}}
\newcommand{\ts}{\tilde{\sigma}}
\newcommand{\kk}{\mathbf{k}}
\newcommand{\tOmega}{\tilde{\Omega}}
\newcommand{\tomega}{\tilde{\omega}}
\newcommand{\hE}{\hat{E}}
\newcommand{\hb}{\hat{b}}
\newcommand{\Primes}{R}
\newcommand{\ol}{\overline{\ell}}
\newcommand{\oM}{{\overline{M}}}
\newcommand{\oMdual}{\overline{M}{}^\dual}
\newcommand{\loc}{\textrm{loc}}
\newcommand{\glob}{\textrm{glob}}
\newcommand{\smooth}{\textrm{smooth}}
\newcommand{\shared}{\textrm{sh}}
\newcommand{\dual}{*}

\begin{document}
\title{Characterizing Generic Global Rigidity}

\author[Gortler]{Steven J. Gortler}
\address{School of Engineering and Applied Sciences, Harvard
  University, Cambridge, MA 02138}
\email{sjg@cs.harvard.edu}

\author[Healy]{Alexander D. Healy}
\address{School of Engineering and Applied Sciences, Harvard
  University, Cambridge, MA 02138}
\email{ahealy@post.harvard.edu}

\author[Thurston]{Dylan P. Thurston}
\address {Department of Mathematics, Barnard College, Columbia University\\ New York, NY 10027}
\email {dthurston@barnard.edu}

\begin{abstract}
  A $d$-dimensional \emph{framework} is a graph and a map from its vertices
  to~$\EE^d$.  Such a framework is \emph{globally rigid} if
  it is the only framework in $\EE^d$ with the same graph and edge
  lengths, up to
  rigid motions.  For which underlying graphs
  is a generic framework globally rigid?  We answer this question
  by proving a conjecture by Connelly, that his sufficient condition
  is also necessary: a generic framework is globally rigid if and only
  if it has a stress matrix with kernel of dimension $d+1$, the
  minimum possible.

  An alternate version of the condition comes from considering the
  geometry of the length-squared mapping~$\ell$: the
  graph is generically locally rigid iff the rank of $\ell$ is maximal,
  and it is generically
  globally rigid iff the rank of the Gauss map on the image of $\ell$
  is maximal.

  We also show that this condition is efficiently
  checkable with a randomized algorithm, and prove that if a graph is
  not generically globally rigid then it is flexible one dimension
  higher.
\end{abstract}

%\keywords{rigidity of graphs, Gauss map, generic frameworks}

\maketitle

%\tableofcontents

\section{Introduction}
\label{sec:intro}

In this paper we characterize those generic frameworks which are
globally rigid in $d$\hyp dimensional Euclidean space. We do this by
proving a conjecture of
Connelly~\cite{Connelly05:GenericGlobalRigidity}, 
who described a 
sufficient condition for generic 
frameworks to be globally
rigid, and conjectured that this condition was necessary.  
As this
condition depends only on the graph and the dimension $d$,
and not on the specific (generic) framework of that graph, we can conclude that
generic global rigidity in $\EE^d$ is a property of a graph.
We further
show that this property
can be checked in probabilistic polynomial
time.

Global rigidity has applications in chemistry, where
various 
technologies measure inter-atomic distances. From this data
one may try
to infer the geometry of the 
configuration \cite[inter alia]{crippen}. This inference
is only well posed if the associated framework is globally
rigid. Moreover, testing for generic global rigidity
can be used as part of a divide-and-conquer strategy for this 
inference~\cite{Hendrickson95:MoleculeProblem}. Similar problems
arise in the field of sensor networks \cite[inter alia]{hamdi}.

\subsection{Definitions and results}

\begin{definition}\label{def:config-space}
  A \emph{graph}~$\Gamma$ is a set of $v$ vertices $\Verts(\Gamma)$ and
  $e$ edges~$\Edges(\Gamma)$, where $\Edges(\Gamma)$ is a set of two\hyp
  element subsets of
  $\Verts(\Gamma)$.
We will typically drop the graph~$\Gamma$ from this notation.
A \emph{configuration} of $\Gamma$ in $\EE^d$
is a mapping from $\Verts(\Gamma)$  to Euclidean space~$\EE^d$.
A \emph{framework}~$\rho$ in $\EE^d$
is a graph~$\Gamma$ together with a configuration of $\Gamma$ in $\EE^d$; we
will also say that $\rho$ is a framework of~$\Gamma$.
Let $C^d(\Gamma)$ denote
the space of frameworks with a given graph~$\Gamma$ and dimension~$d$.
For $\rho\in C^d(\Gamma)$ and  $u \in
\Verts(\Gamma)$,
let $\rho(u)$ denote the image of $u$ under the configuration of~$\rho$.
For a given graph~$\Gamma$ and dimension~$d$,
the \emph{length-squared function} $\ell:C^d(\Gamma)\rightarrow\RR^e$ 
is the  function assigning to each edge of $\Gamma$ its squared 
edge length in the framework.
In particular, the component of $\ell(\rho)$ in the direction of an edge
  $\{u,w\}$ is $\abs{\rho(u)-\rho(w)}^2$.
\end{definition}

\begin{definition}\label{def:globally-rigid}
  A framework~$\rho$ in $\EE^d$ is \emph{congruent} to another
  framework if they are related by an element of the group $\Eucl(d)$
  of rigid motions of $\EE^d$ (rotations, reflections, and
  translations).  We say that $\rho$ is \emph{globally rigid} 
if $\rho$ is the
  only framework of $\Gamma$
in $\EE^d$
with the same edge lengths, up to congruence.
  Equivalently, $\rho$ is
  globally rigid iff 
  $\ell^{-1}(\ell(\rho))/\Eucl(d)$
consists of $[\rho]$, the congruence class of frameworks which contains~$\rho$.
\end{definition}

\begin{definition}\label{def:locally-rigid}
  A framework $\rho \in C^d(\Gamma)$ is \emph{locally rigid} 
if 
there exists a neighborhood~$U$ of $\rho$ 
in $C^d(\Gamma)$
such that $\rho$  is the
  only framework 
in $U$
with the same set of edge lengths, up 
to congruence; equivalently, $[\rho]$ is isolated in
$\ell^{-1}(\ell(\rho))/\Eucl(d)$.%
\end{definition}

(This property is usually just called 
\emph{rigidity} but we will use term \emph{local rigidity} in this paper to make
explicit its distinction from  \emph{global rigidity}.)

\begin{definition}
  \label{def:generic}
  A framework is \emph{generic} if the coordinates of its
configuration do
  not satisfy any non-trivial algebraic equation with rational 
  coefficients.
\end{definition}

We first deal with graphs with very few vertices.
Asimow and Roth proved that a generic framework $\rho$ in $\EE^d$
of a
graph~$\Gamma$ with $d+1$ or fewer vertices
 is  globally rigid 
if $\Gamma$ is a complete graph (i.e.,
a simplex), otherwise it is not even 
locally rigid~\cite[Corollary 4]{AR78:RigidityGraphs}.
Therefore in the rest of the paper we may assume that our graph has $d+2$
or more vertices. In particular, this implies that a
generic framework
does not lie
in a proper affine subspace of $\EE^d$.

Since $\Eucl(d)$ acts freely on frameworks
that do not lie
in a proper affine subspace of~$\EE^d$,
for such a framework~$\rho$, the fiber $\ell^{-1}(\ell(\rho))$
always has
dimension at least $\dim(\Eucl(d))$ (which is $\binom{d+1}{2}$).
In particular, a generic framework with at least $d+1$ vertices
satisfies this condition.
The intuition behind the following definition and
theorem is that the
dimension of a generic fiber of an algebraic map~$f$ is the same as
the kernel of the Jacobian~$df_x$ of~$f$ at a generic point~$x$, and
for a generic locally rigid framework the kernel of $df_x$ only
contains the tangents to the action of $\Eucl(d)$.

\begin{definition}
  Let $d\ell_\rho$ be the \emph{rigidity matrix} of~$\rho$, the
  Jacobian of $\ell$ at the framework $\rho$; by definition, this is a
  linear map from $C^d(\Gamma)$ to $\RR^e$.
  A framework $\rho \in C^d(\Gamma)$ of a graph $\Gamma$ with $d+1$ or
  more vertices is \emph{infinitesimally rigid} if
  \begin{equation}
    \rank d\ell_{\rho} = \dim (C^d(\Gamma)) - \dim (\Eucl(d)) = vd -
    \binom{d+1}{2}.
  \end{equation}
\end{definition}

\begin{theorem}[Asimow-Roth~\cite{AR78:RigidityGraphs}]
  \label{thm:locally-rigid}
  If a generic
framework $\rho \in C^d(\Gamma)$ of a graph $\Gamma$ with $d+1$ or more vertices
is locally rigid
in $\EE^d$, then it is infinitesimally rigid.
Otherwise $\rank d\ell_\rho$ is lower than $vd - \binom{d+1}{2}$.
\end{theorem}

Since the rank of 
the linearization an algebraic map is the same (and maximal) at every
generic point, 
  local rigidity 
in $\EE^d$
is a generic property. 
(See also Lemma~\ref{lem:poly-mat}.)
That is, either all generic frameworks in
$C^d(\Gamma)$ are locally rigid,
or none of them are.
Thus we can call this condition
\emph{generic local rigidity in $\EE^d$} and consider it as 
a property of the  graph $\Gamma$.

We next define some concepts we need to state the condition for
generic global rigidity.

\begin{definition}
  \label{def:stress-vector}
  An 
 \emph{equilibrium stress vector} of a framework~$\rho$ of $\Gamma$ is a 
  real valued function~$\omega$
  on the (undirected) 
  edges of~$\Gamma$ so that 
for all
  vertices $u \in \Verts$,
  \begin{equation}\label{eq:stress-def}
  \sum_{\{w \in \Verts \mid \{u,w\} \in \Edges\}}\!\!\omega(\{u,w\})(\rho(w) -
  \rho(u)) = 0.
  \end{equation}
  (Note that the all-zero function is an equilibrium 
 stress vector for any~$\rho$.)
  Let $S(\rho)$ be the vector space of equilibrium 
stress vectors of the
  framework~$\rho$.
\end{definition}

\begin{remark}
  In many cases, we can interpret \eqref{eq:stress-def} as saying that
  each vertex~$u$ is the weighted average of its neighbors, with
  weights given by $\omega(\{u,w\})$.  This interpretation breaks down
  if $\sum_w \omega(\{u,w\}) = 0$, as happens, for instance, for
  $K_{5,5}$ in~$\EE^3$.  (See Section~\ref{sec:k_5.5-examp}.)  An
  alternate interpretation is that if we put springs on the edges
  of~$\Gamma$ with spring constants given by~$\omega$, the resulting
  forces on each vertex balance out.  Note, however, that some of the
  spring constants in the graph are necessarily negative, and the rest
  position is when the spring is at zero length.
\end{remark}

\begin{definition}
  \label{def:stress-matrix}
  An
 \emph{equilibrium stress matrix} of a framework~$\rho$ of $\Gamma$
  is a rearrangement of an equilibrium 
 stress vector into a
  matrix form, with suitably chosen diagonal entries.
More precisely, it is a matrix~$\Omega$
  indexed by $\Verts\times \Verts$ so that
  \begin{enumerate}
  \item for all $u,w \in \Verts$, $\Omega(u,w) = \Omega(w,u)$;
  \item for all $u,w \in \Verts$ with $u \ne w$ and $\{u,w\} \not\in \Edges$,
    $\Omega(u,w) = 0$;
  \item\label{item:stress-row-sum} for all $u \in \Verts$, $\sum_{w\in\Verts} \Omega(u,w) = 0$; and
  \item\label{item:cond-balance} for all
    $u \in \Verts$,
    $\sum_{w\in\Verts} \Omega(u,w)\rho(w) = 0$.
  \end{enumerate}
  For every equilibrium 
stress vector~$\omega$ there is a unique  equilibrium stress
  matrix~$\Omega$ so that for $u \ne w$ and $\{u,w\} \in \Edges$,
  $\Omega(u,w) =
  \omega(\{u,w\})$, and conversely.
  (This implies that $\Omega(u,u) = -\sum_{w \ne u}
  \omega(\{u,w\})$.)  Thus we will freely convert between equilibrium 
stress
  vectors~$\omega$ and equilibrium stress matrices~$\Omega$. If
  $\Omega$ is an  equilibrium
 stress matrix for $\rho$,
we will also say that $\rho$ \emph{satisfies} the
stress matrix~$\Omega$.
\end{definition}

\begin{definition}
  \label{def:stress-kernel}
  For $\Omega$ an equilibrium stress matrix,
  the \emph{stress kernel}~$K(\Omega)$ is the kernel of $\Omega$.  
  It is isomorphic to the space of frameworks of $\Gamma$ in $\EE^1$ which
  have~$\Omega$ as an equilibrium stress.
  Let $k(\Omega)$ or simply $k$ be $\dim(K(\Omega))$.

Let $k_{\min}(\Gamma, d)$ or just $k_{\min}$ be the minimal value of
$k(\Omega)$ as $\Omega$ ranges over all stress matrices of all
generic frameworks in $C^d(\Gamma)$.  
As has been noted by Hendrickson~\cite[Theorem
2.5]{Hendrickson95:MoleculeProblem}, for each fixed generic framework
$\rho$, the minimal value of $k(\Omega)$ as $\Omega$ ranges over
$S(\rho)$ agrees with this $k_{\min}$.  Indeed 
$k(\Omega) = k_{\min}$ for all generic frameworks~$\rho$ and generic
$\Omega \in S(\rho)$,
suitably defined.
(See
Lemma~\ref{lem:generic-stress-rank}.)
\end{definition}

$K(\Omega)$ always contains the subspace
spanned by the coordinates of~$\rho$ along each axis and the
vector~$\vec 1$ of all $1$'s. This corresponds to the fact that any
affine image of $\rho$ satisfies all of the  equilibrium 
stress matrices in  $S(\rho)$.
If $\rho$ does not lie in a proper affine subspace
of~$\EE^d$, these vectors are independent and span a
$(d+1)$-dimensional space.  We have therefore proved the following
lemma.

\begin{lemma}
  For frameworks of a graph~$\Gamma$ with at least $d+1$
  vertices, $k_{\min}(\Gamma, d) \geq d+1$.
\end{lemma}

\begin{definition}
  \label{def:deficient}
  A graph~$\Gamma$ has 
a \emph{minimal stress kernel in $\EE^d$}
if $k_{\min}(\Gamma, d) = d+1$.
\end{definition}

We can now state Connelly's criterion for 
generic global rigidity.

\begin{theorem}[Connelly~\cite{Connelly05:GenericGlobalRigidity}]
\label{thm:connelly-rigid}
  If a graph $\Gamma$ with $d+2$ or more 
vertices has a minimal stress kernel in $\EE^d$, then 
all generic frameworks $\rho \in C^d(\Gamma)$
are globally rigid. 
\end{theorem}

Our main result is the converse of Connelly's result, which completes the
characterization of generic global rigidity.

\begin{theorem}
  \label{thm:deficient-flexible}
  If a graph $\Gamma$ with $d+2$ or more 
vertices does not have a minimal stress kernel in $\EE^d$, then 
any generic framework $\rho \in C^d(\Gamma)$ is not
globally rigid.
\end{theorem}

The proof of Theorem~\ref{thm:deficient-flexible} is in
Section~\ref{sec:proof}.

Putting these two together, we can also conclude that global rigidity
is a generic property.

\begin{corollary}\label{cor:generic-prop}
Either all generic frameworks in $C^d(\Gamma)$ of a graph are globally
rigid, or none of them are.
Thus we can call this condition {\em
generic global rigidity in $\EE^d$} and consider it as  
a property of a graph.
\end{corollary}

Since the condition is an algebraic condition,
equivalent statements to generic global rigidity in $\EE^d$
include the following:
\begin{itemize}
\item frameworks in a
Zariski open subset of $C^d(\Gamma)$ are globally rigid; or
\item frameworks
in a subset of
full measure of $C^d(\Gamma)$ are globally rigid.
\end{itemize}

In fact, we can deduce a stronger statement.

\begin{corollary}
  If a graph $\Gamma$ with $d+2$ or more vertices does not have
a minimal stress kernel in $\EE^d$, 
then any infinitesimally rigid
  framework $\rho \in C^d(\Gamma)$ is not globally rigid.%
\end{corollary}

(This corollary can also be proven by combining
corollary~\ref{cor:generic-prop} with
\cite[Proposition 2.4]{CW05:TransferGlobalRigidity}.)

\begin{proof}
  Let $\rho_0$ be an infinitesimally rigid framework.  First note that
  $\rho_0$ does not lie in any affine subspace of $\EE^d$, so by
  Proposition~\ref{prop:sing-codim-Arho} it maps
  to a smooth point $[\rho_0]$ in the quotient $C^d(\Gamma)/\Eucl(d)$.
  By the inverse function theorem, $[\rho_0]$ has
  an open neighborhood~$N$ that maps
  injectively into the space of edge lengths.  Consider a
  sequence~$\{[\rho_i]\}$ of classes of generic configurations inside $N$
  approaching~$[\rho_0]$.  By Theorem~\ref{thm:deficient-flexible}, for
  each $\rho_i$ there is an inequivalent configuration $\sigma_i$ with
  the same edge lengths.  Because the map $\ell$ is proper on
  $C^d(\Gamma)/\Eucl(d)$ (see
  Lemma~\ref{lem:lisproper}), the
  sequence $\{[\sigma_i]\}$ has an accumulation point~$[\sigma_0]$, with
  $\ell(\sigma_0) = \ell(\rho_0)$ by continuity.  Since
  $\ell$ is injective on $N$, we have  $[\sigma_i] \not\in N$ and so
  $[\sigma_0] \not\in N$; thus $\sigma_0$ is the desired
  framework with the same edge lengths as $\rho_0$ but not congruent
  to~$\rho_0$. 
\end{proof}

Note that there is no hope for a similar strengthening of
Theorem~\ref{thm:connelly-rigid}; global rigidity can fail at
non-generic points in any number of ways, for instance if there is a
transverse self-intersection of the image space $\ell(C^d(\Gamma))$.

As a corollary of Theorems~\ref{thm:connelly-rigid}
and~\ref{thm:deficient-flexible}, we find an efficient algorithm for
checking for generic global rigidity.

\begin{theorem}\label{thm:random-global-RP-1}
  There is a polynomial-time randomized algorithm for checking for generic
  global rigidity in~$\EE^d$.
\end{theorem}

Essentially, 
you take random a 
framework~$\rho$ in $C^d(\Gamma)$,
using sufficiently large random numbers as coordinates,
and a
random equilibrium stress matrix~$\Omega$ of~$\rho$. Then with high
probability, $\dim (\ker(\Omega)) = k_{\min}$.
In Section~\ref{sec:complexity} we prove a more precise version of
this statement (Theorem~\ref{thm:random-global-RP}), in particular
giving bounds on ``sufficiently large''.

One further result, arising from similar considerations involved in
the proof of Theorem~\ref{thm:deficient-flexible}, is the
following:

\begin{theorem}\label{thm:flex}
  If a graph~$\Gamma$ is not generically
  globally rigid in~$\EE^d$, then any generic framework~$\rho$ in
  $C^d(\Gamma)$ can be connected to some incongruent framework~$\sigma$ 
  in $C^d(\Gamma)$ by a
  path of frameworks of $\Gamma$ in $\EE^{d+1}$ with constant edge lengths.
\end{theorem}

Theorem~\ref{thm:flex} is proved in Section~\ref{sec:flexes}.

\subsection{Relationship to previous work}

The problem of $d$-embeddability, determining if there is a
framework of a graph
in $\EE^d$ with a given set of (integer) edge lengths, was shown to be strongly
NP-hard by Saxe~\cite{Saxe79:EmbedGraphsNP}. Saxe also showed that the
problem of determining 
global rigidity is strongly NP-hard.  The proof starts with
the observation that a 1-dimensional framework of a graph formed of
$n$~vertices connected in a cycle is a solution to the partition
problem on the edge lengths.

But there has been some hope that global rigidity can be characterized
and efficiently tested if one suitably restricts the input to avoid certain
coincidences (like the coincidence that there be two solutions to a
partition problem that has at least one solution).
In particular, one strong way of restricting the
problem is to assume that the original framework is generic.  For
example, it is
easy to see that a generic framework of a graph
is globally rigid in $\EE^1$ iff 
the graph 
is 2-connected.

Hendrickson~\cite{Hendrickson95:MoleculeProblem} describes the 
condition of Theorem~\ref{thm:connelly-rigid}
ascribed to Connelly (CC for short), for generic frameworks of
a graph to be globally rigid in
$\EE^d$. Connelly proved the sufficiency of this
condition~\cite{Connelly05:GenericGlobalRigidity}.
Connelly also conjectured that CC was a necessary
condition, our Theorem~\ref{thm:deficient-flexible}.

Hendrickson also gave~\cite{Hendrickson92:ConditionsUniqueGraph} two
necessary conditions (HC for short), for
generic frameworks of~$\Gamma$ to be globally rigid in $\EE^d$.
Namely, $\Gamma$ must be $(d+1)$-connected
and deleting any edge of $\Gamma$ must leave a generically locally rigid graph.
A combinatorial
result of Jackson and Jord\'an~\cite{JJ05:ConnRigMatriodGraphs}
effectively shows that for $d=2$, $\text{CC}=\text{HC}$, and
so both conditions are necessary and sufficient.  They later gave
a direct proof of the sufficiency of
HC for $d=2$~\cite{JJ06:GloballyLinked}.  Moreover in 2 dimensions, 
generic local
rigidity and therefore generic global rigidity can be characterized
combinatorially~\cite{laman} and
deterministically solved for efficiently \cite[inter alia]{lovasz}.

Connelly showed~\cite{Connelly91:GenericGlobalRigidity} that HC are
not sufficient to show that  a generic framework is 
globally rigid  for $d \ge 3$; specifically, he showed that the complete
bipartite graph
$K_{5,5}$, despite satisfying
HC,  is not globally rigid in 3 dimensions for some
generic frameworks.  Similar results hold for larger bipartite graphs
in higher dimensions.
These appear to be the only known examples of such graphs.

(See Sections~\ref{sec:not-redund-rigid-examp}
and~\ref{sec:d-conn-examp} for more on graphs that fail HC,
and Section~\ref{sec:k_5.5-examp} for more on $K_{5,5}$.)

\subsection{Structure of the paper}
We prove the main result in Section~\ref{sec:proof} and give a
few examples of the relevant spaces in Section~\ref{sec:examples}.

In Section~\ref{sec:shared-stress} we give an alternate
characterization of generic global rigidity in terms of the
\emph{shared stress kernel} (the intersection of the kernels of all
stress matrices) or, alternatively, the rank of the \emph{Gauss map}.
(The Gauss map is the map that takes a smooth point of an algebraic
variety to its tangent space, considered as a point in an appropriate
Grassmannian.)  We also relate the theorem to the more general
algebraic geometry setting of secant varieties and Zak's theorem.

We then show that
the condition of generic global rigidity can be
checked in randomized polynomial time in
Section~\ref{sec:complexity}, following remarks of
Hendrickson~\cite{Hendrickson92:ConditionsUniqueGraph,Hendrickson95:MoleculeProblem}.

In Section~\ref{sec:flexes} we prove Theorem~\ref{thm:flex}.

Throughout this paper, we attempt to give enough details for readers
who are not very familiar with algebraic geometry, while at the same
time indicating how the more abstract arguments work.

\subsection*{Acknowledgements}
We would like to thank Robert Connelly, Daniel Freedman, Yael Karshon,
J. M. Landsberg, Amit Singer, Ileana Streinu and an anonymous referee for helpful
conversations and suggestions.
AH was supported by NSF grant CCR-0205423 and a Sandia fellowship.
DPT was supported by a Sloan research fellowship.

%%% Local Variables: 
%%% mode: latex
%%% TeX-master: "main"
%%% End: 

\section{Proof of main theorem}
\label{sec:proof}

Our basic strategy for the proof of
Theorem~\ref{thm:deficient-flexible} is similar in spirit to the
approach used by
Hendrickson~\cite{Hendrickson92:ConditionsUniqueGraph} to show the
necessity of what he called ``redundant rigidity'' for  a generic 
global rigidity in $\EE^d$. Given a graph~$\Gamma$ which does not have
a minimal stress kernel in~$\EE^d$,
we
construct spaces $X$ and~$Y$ and a map $f:X \rightarrow Y$ so that the
preimages in $X$ of a point in $Y$ correspond to incongruent
frameworks 
with the same edge lengths.
We then show that the
``degree mod-two'' of this map is well defined and equal to~0. This
degree is equal to the number of preimages at a regular value, modulo
2.  Thus for
such a
map, any regular value in the image must have at least one more
preimage, which represents a distinct framework 
in $C^d(\Gamma)$
with the same edge
lengths, thus contradicting global rigidity in $\EE^d$.  To guarantee a well
defined degree, we will not have the
luxury of the domain being a compact manifold. Rather the domain will be
a smooth space with singularities of co-dimension $2$ or greater,
and $f$ will be a \emph{proper map}.  (See
Definition~\ref{def:proper} and
Theorem~\ref{thm:degree}.)

The main spaces involved in the proof are summarized
in Figure~\ref{fig:spaces}, as we now explain.
\begin{figure}
  \begin{center}
\begin{tikzpicture}
  \matrix(spaces) [matrix of math nodes, column sep=1cm, row sep=1cm]
  {
    A(\Omega)    & A(\Omega)/\Eucl(d)& B(\Omega)  &  L(\Omega)\\
    C^d(\Gamma)& C^d(\Gamma)/\Eucl(d)& M(\Gamma)& \RR^e\\
  };
  \draw[->>] (spaces-1-1) -- (spaces-1-2);
   \draw[->>] (spaces-1-2) -- node[cdarrow]{\text{same}} node[cdarrow,below]{\text{dim}} (spaces-1-3);
    \draw[right hook->] (spaces-1-3) -- node[cdarrow]{\text{same}} node[cdarrow,below]{\text{dim}} (spaces-1-4);
  \draw[->>] (spaces-2-1) -- (spaces-2-2);
   \draw[->>] (spaces-2-2) -- node[cdarrow]{\text{same}} node[cdarrow,below]{\text{dim}} (spaces-2-3);
    \draw[right hook->] (spaces-2-3) -- (spaces-2-4);
  \draw[right hook->] (spaces-1-1)--(spaces-2-1);
  \draw[right hook->] (spaces-1-2)--(spaces-2-2);
  \draw[right hook->] (spaces-1-3)--(spaces-2-3);
  \draw[right hook->] (spaces-1-4)--(spaces-2-4);
\end{tikzpicture}
  \end{center}
  \caption{The main spaces involved in the proof of Theorem~\ref{thm:deficient-flexible}}
  \label{fig:spaces}
\end{figure}
We start with the map $\ell$ from $C^d(\Gamma)$ to $\RR^e$, and
factor it through the quotient by the Euclidean group, followed by the
surjection onto its image $M(\Gamma) \subset \RR^e$.  The spaces
$C^d(\Gamma)/\Eucl(d)$ and its image $M(\Gamma)$ have the same
dimension by local rigidity.

We next take any generic~$\rho$ and pick an equilibrium stress~$\Omega$
for~$\rho$ which is
generic in a suitable sense, and consider the space $A(\Omega)$ of all
frameworks satisfying~$\Omega$.  As motivation for this, recall that
Connelly's proof of Theorem~\ref{thm:connelly-rigid} first shows that
for any generic framework~$\rho$ and a global flex~$\sigma$ of it, the
map~$\ell$ agrees on neighborhoods of~$\rho$ and~$\sigma$
\cite[Proposition 3.3]{Connelly05:GenericGlobalRigidity}.  In
particular, the tangent space to~$M(\Gamma)$ at~$\rho$ and~$\sigma$
agree, so any equilibrium stress for~$\rho$ is also an
equilibrium  stress for~$\sigma$.
(Compare Lemma~\ref{lem:tangent-stress}.)
Thus $\sigma$ must lie in $A(\Omega)$, for any $\Omega \in S(\rho)$.
So we may as well
look for a global flex within such a space $A(\Omega)$.

We  define  the  domain  of  our map  to  be  $A(\Omega)/\Eucl(d)$.   In
Proposition~\ref{prop:sing-codim-Arho},  we show  that  if the  stress
kernel  of $\Omega$
is of  dimension  greater  than $d+1$,  then  this domain  has
singularities  of   codimension  at  least~$2$.   This   is  based  on
recognizing $A(\Omega)$  as $K(\Omega)^d$,  where $K(\Omega)$ is the stress
kernel as in Definition~\ref{def:stress-kernel}.

We also consider the image of $A(\Omega)$ under~$\ell$, denoted
$B(\Omega) \subset M(\Gamma)$.
Again, by
local rigidity, the quotient $A(\Omega)/\Eucl(d)$ and $B(\Omega)$ have
the same dimension.
If we view the equilibrium stress vector as defining a hyperplane in $\RR^e$,
we can then alternately view $B(\Omega)$ as the \emph{contact locus}
of the corresponding equilibrium stress~$\omega$: the set of points in
$M(\Gamma)$ whose tangent space is
contained in the hyperplane defined by~$\omega$.  A fundamental
theorem of algebraic geometry says that a generic contact locus is
\emph{flat}: it is contained in a linear space $L(\Omega)$ of the same
dimension.  (See Proposition~\ref{prop:contact-flat}.)  The space $L(\Omega)$ serves
as the range.

Our map $\ol$ is the restriction of $\ell$ to this domain and range.
Since the singularities of the domain have
codimension at least~$2$, $\ol$ will have a well-defined degree.
The final step in the proof is that this degree is~$0$, which is
guaranteed by noticing that any point in $\RR^e$ with at least one
negative coordinate is not in the image of~$\ell$.

\medskip

If a generic framework in $C^d(\Gamma)$
is not locally rigid,
then it is clearly not  globally rigid, in which case
our proof is done.
So, in what follows, we will assume that the framework $\rho\in C^d(\Gamma)$
is locally rigid
and has at least $d+2$ vertices, and so by
Theorem~\ref{thm:locally-rigid},
$\rank d\ell_{\rho} = vd - \binom{d+1}{2}$.

\subsection{Algebraic geometry preliminaries}
\label{sec:algebr-geom}

We start with some preliminaries about semi-algebraic sets from real
algebraic geometry, somewhat specialized to our particular case.  For
a general reference, see, for instance, the
book by Bochnak, Coste, and Roy~\cite{bcr}.
\begin{definition}
  \label{def:generic-gen}
  An affine, real \emph{algebraic set} or \emph{variety}~$V$ 
  defined over a field~$\kk$
  contained in~$\RR$ is a subset of $\RR^n$ that is
  defined by a set of algebraic equations with coefficients in~$\kk$.
  It is \emph{irreducible} if it is not the union of two proper
  algebraic subsets
  defined over $\RR$. 
  An algebraic set has a \emph{dimension}
  $\dim(V)$, which we will define as the largest $t$ for which there
  is an open (Euclidean) subset of~$V$ isomorphic to $\RR^t$.

  A point~$x$ of an irreducible algebraic set~$V$ is \emph{smooth} if it has a
  neighborhood that is
  smoothly isomorphic to $\RR^{\dim(V)}$; in this case there is a
  well-defined tangent space $T_x(V)$.  (Note that there may be points
  with neighborhoods isomorphic to $\RR^n$ for some $n < \dim(V)$; we
  do not consider these points to be smooth.)
\end{definition}

\begin{definition}
  A \emph{semi-algebraic set}~$S$ over $\kk$ 
  is a subset of $\RR^n$ defined by algebraic
  equalities and inequalities with coefficients in $\kk$; 
  alternatively, it is the image of an
  algebraic set (defined only by equalities) under an algebraic map
  with coefficients in $\kk$.
  A semi-algebraic set has a well defined (maximal) dimension~$t$.

  The \emph{Zariski closure} of $S$ is the smallest
  algebraic set over $\RR$ containing it.  (Loosely speaking, we can get an 
algebraic set by keeping all algebraic
  equalities and dropping  the inequalities. We may 
need to enlarge the  field to cut out  the smallest algebraic set
containing~$S$
but a finite extension will always suffice.)
  We call $S$ \emph{irreducible} if its Zariski closure is
  irreducible. (We chose here to define the Zariski closure as well as the
notion of irreducibility using the field $\RR$ instead of $\QQ$ or
$\kk$ in order
to avoid
complications in the proof of Proposition~\ref{lem:conormal-props} below.)
  An irreducible semi-algebraic set $S$ has the same dimension as its
  Zariski closure.  

A point on $S$ is smooth
  if it has a neighborhood in $S$ smoothly isomorphic to
  $\RR^{\dim(S)}$. We denote by $\smooth(S)$ the locus
of smooth points of $S$. It too is semi-algebraic.
We call $S$ itself
smooth if $\smooth(S)=S$.

\end{definition}

Suppose that $S$ is a smooth irreducible semi-algebraic set.
Let~$\phi$ be a non-zero algebraic function on~$S$.  Its
zero set is stratified as a union of a finite number of smooth 
manifolds each with dimension less than 
$\dim(S)$, so (since $S$ is smooth) the complement of the zero set
(i.e., $\{\,x\in S \mid\phi(x) \ne 0\,\}$) is open and dense in $S$.

\begin{definition}
  The image of~$\ell$ is called the \emph{measurement set} and denote
  by~$M(\Gamma)$ or just~$M$.
  It is a  semi-algebraic set defined over~$\QQ$.
\end{definition}

\begin{lemma}
  \label{lem:M-irred} The semi-algebraic set~$M$ is irreducible.
\end{lemma}

\begin{proof}
  $M$ is the image of $C^d(\Gamma)$ under a polynomial
  map.  Since $C^d(\Gamma)$ is irreducible and the image of an
  irreducible variety is irreducible, it
  follows that~$M$ is irreducible.
\end{proof}

The dimension of $M$ is the rank of $d\ell_{\rho}$ at any
generic configuration $\rho$ (and is $vd - \binom{d+1}{2}$ when the
graph is generically locally rigid in $\EE^d$).  More precisely, the
tangent space to $M$ for generic~$\rho$ is the span of the
image of $d\ell_{\rho}$.

We next define genericity in larger generality and give some basic
properties.

\begin{definition}
  A point in
  a  (semi-)algebraic set~$V$ defined over~$\kk$ is
  \emph{generic} if its
  coordinates do not satisfy
  any algebraic equation with coefficients in~$\kk$
  besides those that are satisfied by every
  point on~$V$.

  Almost every point in an irreducible semi-algebraic set $V$ is generic.
\end{definition}

Note that the defining field might change when we take the Zariski
closure.  However, this does not matter for the purposes of genericity.
More formally, if $\kk$ is a finite algebraic extension of $\QQ$
and $x$ is a generic point in an irreducible semi-algebraic set $S$ defined
over~$\kk$, 
then $x$
  is also generic in the Zariski closure of $S$ over an appropriate field.
This follows from a two step argument. 
First, an elementary argument from Galois theory proves that the
(real) Zariski
closure~$V$ of $S$ can be defined using polynomials over a field~$\kk'$
which is some finite extension of~$\kk$.
Second, another elementary argument from Galois theory proves that 
if a point $x$
satisfies
  an algebraic equation with coefficients in~$\kk'$
that is non-zero over $V$,
then $x$ must also 
satisfy
  some algebraic equation with coefficients in~$\kk$
that is non-zero over~$V$.

The following proposition is standard:

\begin{proposition}\label{prop:generic-smooth}
  Every generic point of a (semi-)algebraic set is smooth.
\end{proposition}

\begin{lemma}
  \label{lem:image-generic}
  Let $V\!$ and $W\!$ be irreducible semi-algebraic sets and $f : V \to W$
  be a surjective algebraic map, all defined over~$\kk$.  Then if $x_0
  \in V\!$ is generic, $f(x_0)$ is generic inside~$W\!$.
\end{lemma}

\begin{proof}
  Consider any non-zero algebraic function~$\phi$ on~$W$
  defined over~$\kk$.  Then $x
  \mapsto \phi(f(x))$ is a function on~$V$ that is not identically
  zero.  Thus if $x_0$ is a generic point in $V$, $\phi(f(x_0)) \ne
  0$.  Since this is true for all~$\phi$, it follows that $f(x_0)$ is
  generic.
\end{proof}

Thus, instead of speaking of a generic configuration in $C^d(\Gamma)$,
we may speak of a generic point on $M$: if a configuration~$\rho$ is
generic in $C^d(\Gamma)$, then $\ell(\rho)$ is generic in~$M$.

\begin{lemma}
  \label{lem:preimage-generic}
Let $\kk$ be a finite algebraic extension of $\QQ$.
  Let $V\!$ and $W\!$ be irreducible semi-algebraic sets with $V\!$ smooth,
  and let $f : V \to W$
  be a surjective algebraic map, all defined over~$\kk$.
  Then
  if $y_0\in W$ is generic, there is a point in $f^{-1}(y_0)$ that
  is generic in~$V\!$.
\end{lemma}

\begin{proof}
  Let~$\phi$ be a non-zero algebraic function on~$V$ defined over $\kk$.
  We start by showing there is a point $x \in f^{-1}(y_0)$ so that
  $\phi(x) \ne 0$.  
Consider the semi-algebraic set $X_\phi \coloneqq
  \{\,x \in V \mid \phi(x) \ne 0\,\}$.  This is dense in~$V$ due
to irreducibility and smoothness,
so its
  image $f(X_\phi)$ is dense in~$W\!$.  Therefore $Y_\phi \coloneqq W
  \setminus f(X_\phi)$ is a semi-algebraic subset of~$W$ that is
  nowhere dense.  It follows that there is a non-zero algebraic
  function~$\psi$ on~$W$ defined over~$\kk$ that vanishes when
  restricted to $Y_\phi$.  But then since $y_0$ is generic, $\psi(y_0)
  \ne 0$, which means
  that $y_0$ is in the image of $X_\phi$, so there is an $x \in
  f^{-1}(y_0)$ so that $\phi(x) \ne 0$, as desired.

  Let $Z_\phi = \{\,x\in f^{-1}(y_0) \mid \phi(x) = 0\,\}$.  We have
  shown $Z_\phi$ is a proper subset of $f^{-1}(y_0)$ for all non-zero
  algebraic functions~$\phi$ on $V$.  It follows that for any finite
  collection of $\phi_i$, the union of the $Z_{\phi_i}$ is still a
  proper subset of $f^{-1}(y_0)$ (as we can consider the
  product of the $\phi_i$).  But there are only countably many
  possible $\phi$ overall, and a countable union of algebraic subsets
  covers an algebraic set iff some finite collection of them do.
  (Proof: this is true for each irreducible component, as a proper
  algebraic subset has measure zero, and there are only finitely many
  irreducible components.)  Thus
  the union of the $Z_\phi$ do not cover $f^{-1}(y_0)$, i.e., there is
  a generic point in $f^{-1}(y_0)$.
\end{proof}

\begin{lemma}
  \label{lem:generic-dense}
  Let $V \subset W$ be an inclusion of real algebraic varieties, where
  $W$ (but not~$V$) is defined over~$\QQ$.  Suppose that $V$ is 
smooth and 
  irreducible and that it has one
  point~$y$ which is generic in~$W$ (over~$\QQ$).  Then the points in
  $V$ which are generic in $W$ are dense in~$V$.
\end{lemma}

Note that in this lemma the points we consider are \emph{different}
from the points that are generic in~$V$, since $V$ is not defined
using only the rationals.  (If $V$ were defined using the rationals
and were properly contained in~$W$, then it would have no points which
were generic in~$W$: any rational equation defining~$V$ would by
definition vanish for all $x\in V$.)

\begin{proof}
  Let~$\phi$ be a non-zero algebraic function on~$W$ defined over $\QQ$.
  Consider the semi-algebraic set $X_\phi \coloneqq
  \{\,x \in V \mid \phi(x) \ne 0\,\}$.  
This is non empty due to our assumption about the point $y$.
Thus $X_\phi$ is open and dense in~$V\!$, due to 
irreducibility and smoothness.
The set of points in
$V$ which are generic in $W$ is defined as the intersection of the 
$X_\phi$ as $\phi$ ranges over the countable set of possible $\phi$.
Since $V$ is a Baire space, such a countable intersection 
of open and dense subsets must itself 
be a dense subset.
\end{proof}

\subsection{The domain}

To construct the domain,
we start by considering \emph{stress satisfying} frameworks.

\begin{definition}
  Let~$\Omega$ be 
an equilibrium stress matrix of some framework in $C^d(\Gamma)$.
  The space of \emph{stress satisfiers}~$A(\Omega)$ is the space of all
  $d$-dimensional frameworks that satisfy $\Omega$:
  \begin{equation}
  A(\Omega) \coloneqq \{\,\sigma \in C^d(\Gamma) \mid \Omega \in S(\sigma)\,\}.
  \end{equation}
\end{definition}

\begin{lemma}
  \label{lem:A-K-rel}
  The space $A(\Omega)$ is isomorphic to $K(\Omega)^d$ and has
  dimension $k(\Omega)d$.  In particular, if $\Omega \in S(\rho)$,
  then $A(\Omega)$ contains the space $\af(\rho)$ of affine images
  of~$\rho$.
\end{lemma}

\begin{proof}
  Since Equation~\eqref{eq:stress-def} is an equation for each
  coordinate separately, a framework satisfies an
  equilibrium stress iff each of its projections on the coordinate
  axes does.  But $K(\Omega)$ is isomorphic to the space of
  $1$-dimensional frameworks of~$\Gamma$ that have $\Omega$ as an
  equilibrium stress, so it follows that $A(\Omega)$ is
  $K(\Omega)^d$.  The remaining statements are immediate.
\end{proof}

We now study the quotient $A(\Omega)/\Eucl(d)$, which,
for an appropriate choice of $\Omega$, will be the  domain of
our main map.
If $\Gamma$ is generically locally rigid and
$A(\Omega)$ includes
some framework generic in $C^d(\Gamma)$,
this space is of dimension $k(\Omega)d-\binom{d+1}{2}$.
We need to understand the singularities of this quotient, and in
particular show that they have codimension at least~$2$ in appropriate
cases (so that we can apply Corollary~\ref{cor:degree-sing}).
We will reduce the question to the following lemma.

\begin{lemma}
  \label{lem:sing-codim-gen}
  If $n > d$, the quotient space $(\EE^d)^n/\Eucl(d)$ (i.e.,
  $n$\hyp tuples of points in $\EE^d$, quotiented by $\Eucl(d)$ acting
  diagonally) is a smooth stratified
  space with singularities of codimension $n-d$ and higher.
  Furthermore, the
  singularities occur at classes of points in~$\EE^d$ that lie inside
  a proper affine subspace of~$\EE^d$.
\end{lemma}

Here a smooth stratified space is, loosely speaking, a space which is
decomposed into smooth manifolds of differing dimensions, limiting
onto each other in a nice way.  For instance, any semi-algebraic set,
with any type of singularities, has the structure of a smooth
stratified space. For a formal definition and 
discussion, see~\cite{Pflaum02:AnalGeomStratified}.

\begin{proof}
  In general, the quotient of a variety by a group acting properly is
  stratified by
  conjugacy classes of
  stabilizers \cite[Corollary 4.3.11]{Pflaum02:AnalGeomStratified}.
  In particular the points with trivial stabilizer form a manifold.
  Let us therefore identify the points in $(\EE^d)^n$ where the stabilizer
  is non-trivial.
The key observation is that an $n$\hyp tuple is stabilized by a
non-trivial Euclidean transform if and
only if the dimension of its affine span is less than $d$. 
For example, if all points lie in a hyperplane of dimension
$d-1$, then they are stabilized by the  reflection across that plane.
In general, if the affine span of the $n$\hyp tuple is
  codimension~$r$, the stabilizer is isomorphic to the orthogonal
  group in $r$ dimensions, $O(r)$.

The lowest
  codimension case is when $r=1$, in which case the stabilizer is
  $O(1) \isom \ZZ/2$.
  To compute the codimension of these
  $n$\hyp tuples with $r=1$ in $(\EE^d)^n$, pick $d$ out of the $n$ points
  that span the subspace; these points are unrestricted.  The
  remaining $n-d$ points must lie in that subspace, providing one
  constraint for each point, so these $n$\hyp tuples have codimension $n-d$.
  In general, a similar count shows that $n$\hyp tuples with a
  stabilizer of $O(r)$ are of codimension $r(n-d+r-1)$ in $(\EE^d)^n$.

  Inside the quotient $(\EE^d)^n/\Eucl(d)$, the codimension of the singular
  set is decreased by the
  dimension of the stabilizer.  Thus points with stabilizer $O(1)$
  still have codimension $n-d$ in the quotient, and in general the
  the codimension in the quotient of points with stabilizer $O(r)$ is
  \begin{equation}
  r(n-d+r-1) - \binom{r}{2} = r(n-d) + \binom{r}{2} \le n-d.
  \qedhere
  \end{equation}
\end{proof}

\begin{proposition}
  \label{prop:sing-codim-Arho}
Suppose that $\Omega$ is equilibrium stress matrix so that
$k > d+1$.
 Then $A(\Omega)/\Eucl(d)$ is a smooth stratified space with
  singularities of codimension at least 2.
  Furthermore, the singularities occur at classes of frameworks~$[\sigma]$
  where $\sigma$ lies in a proper affine subspace of~$\EE^d$.
\end{proposition}

\begin{proof}
  Remember that $A(\Omega) \isom K(\Omega)^d$ by Lemma~\ref{lem:A-K-rel}.
  We can view this space as $K(\Omega) \tensor
  \RR^d$ with $\Eucl(d)$ acting on the $d$\hyp
  dimensional coordinates.  We can turn this tensor product around:
  instead of thinking of $A(\Omega)$ as a $d$\hyp
  tuple of points in $K(\Omega)$, think of it as a $k$\hyp tuple of points
  in $\EE^d$, with the diagonal action of $\Eucl(d)$.  Concretely, you
  can always find $k$ vertices of the graph whose position determines
  the position of all other vertices of a framework inside
  $A(\Omega)$.  Then $A(\Omega)$ is isomorphic to $(\RR^d)^k$, thought
  of as the positions of these $k$ vertices.
  
  Now by Lemma~\ref{lem:sing-codim-gen}, $A(\Omega)/\Eucl(d)$ is a
  smooth stratified space with singularities of codimension at least
  $k-d$, at frameworks which lie in a proper affine subspace.
  By assumption, $k > d+1$,
  and so this codimension is at least~$2$.
\end{proof}

\subsection{The range}
\label{sec:range}
We now turn to defining the range.
\begin{definition}\label{def:BL}
  Let~$\Omega$ be
  an equilibrium stress matrix of some framework in $C^d(\Gamma)$.
  The space $B(\Omega)$ is defined to be $\ell(A(\Omega))$, the
  image in $\RR^e$ of squared edge lengths of elements in $A(\Omega)$.
  Let $L(\Omega)$ to be the smallest linear subspace of $\RR^e$
  that contains $B(\Omega)$.

\end{definition}
Note that $B(\Omega)$ is a
semi-algebraic
set.  If the graph is generically locally rigid in $\EE^d$
and $A(\Omega)$ contains a framework generic in $C^d(\Gamma)$,
then 
$B(\Omega)$,
like $A(\Omega)/\Eucl(d)$, is of dimension $kd-\binom{d+1}{2}$ (as the
map $\ell$ is locally one-to-one from the quotient).
The space $L(\Omega)$, for an appropriately chosen $\Omega$,
will be the range space for our map. 
Our  main task is to prove that 
for any generic $\rho$ and appropriately chosen 
$\Omega \in S(\rho)$,
$L(\Omega)$ is in fact the same dimension
as $B(\Omega)$, or that $B(\Omega)$ is flat as in the following definition.

\begin{definition}\label{def:flat}
An irreducible semi-algebraic set~$S$ in $\RR^n$ is \emph{flat} if it
is contained in a linear subspace of $\RR^n$
of the same dimension as~$S$.
\end{definition}

Note that this usage of \emph{flat} is
unrelated to the notion of flat families in algebraic geometry, and an
algebraic variety is
flat by this usage iff it is itself a linear space.

The key point in showing that $B(\Omega)$ is flat is
Proposition~\ref{prop:contact-flat}, a standard theorem asserting
that for any algebraic variety~$V$, the contact locus of a
generic hyperplane in the dual variety $V^{\dual}$ is flat.  We recall
these definitions and theorems.
See \cite{FP01:RuledVarieties} for a 
more complete treatment.

\begin{definition}
  We say that a hyperplane $H$ is \emph{tangent} to a homogeneous
  algebraic set $V \subset \RR^n$ at a smooth point $x\in V$ if
  $T_x(V) \subset H$.  For $\phi$ a functional in the dual space
  $(\RR^n)^*$, we say that $\phi$ is tangent to $V$ at~$x$ when
  $T_x(V) \subset \ker(\phi)$.

  For a homogeneous algebraic set~$V$ in $\RR^n$,
  define the \emph{dual variety} $V^{\dual} \subset (\RR^n)^*$ to be
  the Zariski closure of the set of all functionals that are tangent
  to $V$ at a smooth point.  The \emph{conormal bundle} $C_V \subset
  \RR^n \times (\RR^n)^*$ is the Zariski closure of the pairs $(x, \phi)$
  where $x$ is a smooth point of~$V$ and $\phi$ is tangent to $V$
  at~$x$.  There are two natural projections, $\pi_1: C_V \to V$ and
  $\pi_2: C_V \to V^\dual$.  Above the smooth points of~$V$, $C_V$ is
  smooth, although the projection onto $V^*$ need not be.

  For $\phi \in V^{\dual}$, define the \emph{contact locus}
  $V_\phi$ of~$\phi$ to be $\pi_1(\pi_2^{-1}(\phi))$.
  In particular, this contains the smooth points $x \in V$ so
  that $\phi$ is tangent to $V$ at~$x$.
\end{definition}

\begin{remark}
  In Lemma~\ref{lem:tangent-stress} we
  relate the  space of equilibrium stresses with the 
  dual space to the measurement set~$M$.
  In particular, this implies that the
  dual $\oMdual$ is the Zariski closure of $S(\rho)$ as
  $\rho$ ranges over generic frameworks.
  In Proposition~\ref{prop:Bflat} we further identify $L(\Omega)$, for
  appropriate $\Omega$, as a contact locus.
\end{remark}

\begin{lemma}\label{lem:conormal-props}
  For $(x,\phi) \in C_V$ with $x$ smooth in~$V$, 
 we have that $\phi$ is tangent to~$x$.  (That is, the
  Zariski closure in the definition of $C_V$ does not
  add points above smooth points.)  Similarly if  $x \in
  V_\phi$ with $x$ smooth in $V$, $\phi$ is tangent to~$x$.
If $V$ is a homogeneous variety in
  $\RR^n$, the conormal bundle $C_V$ always has
  dimension~$n$ (independent of the dimension of~$V$).  If $V$ is
  irreducible, then so are $C_V$ and $V^*$.
\end{lemma}

\begin{proof}
  See, e.g., \cite[Section 2.1.4]{FP01:RuledVarieties} for the proofs
  of these basic properties in the complex
  case.  For the passage to the real case, let $V^\CC$ be the
  complexification of~$V$ and let $F: \RR^n \to \RR^m$ be an algebraic
  map that cuts out $V$, in the sense that $V$ is $F^{-1}(0)$ and
  $\ker dF_x = T_x(V)$ at all smooth
  points $x \in V$.  Define
  \begin{equation}
    \Gamma_V \coloneqq
    \Bigl\{\,(x,\phi) \in V  \times (\RR^n)^* \mathrel{\Big\vert}
      \rank \begin{pmatrix} \phi \\ dF_x \end{pmatrix}
      \le n - \dim V\,\Bigr\}.
  \end{equation}
  For any smooth point $x \in V$, the set of $\phi$ such that 
  $(x,\phi)$ is in $\Gamma_V$
  corresponds
  exactly to the set of functionals that are tangent to $V$ at $x$.
  Thus $C_V$ must be a union of 
  irreducible components of $\Gamma_V$.

  Most of the properties in the real case follow immediately from
  the observation that 
$(\Gamma_V)^\CC =
  \Gamma_{V^\CC}$, with the possible exception of irreducibility.  For
  irreducibility of $C_V$, 
  notice that 
from the properties of 
a complex conormal bundle, e.g., \cite[Section 2.1.4]{FP01:RuledVarieties},
we know that 
$C_{V^\CC}$ is an irreducible component of
  $\Gamma_{V^\CC}$
  but by, e.g., \cite[Lemma
  7]{Whitney57:ElemStructRealVarieties}, the components of
  $(\Gamma_V)^\CC$ correspond to the components of $\Gamma_V$, so in
  particular there is one component containing all the tangents to
  smooth points.  Irreducibility of $V^*$ follows from projection
  and the fact that from their definitions $V^*$ is the same as the 
  Zariski closure of $\pi_2(C_V)$.
\end{proof}

\begin{proposition}\label{prop:double-dual}
  For a homogeneous irreducible algebraic set~$V$, the double dual
  $V^{\dual\dual}$ is $V$.
\end{proposition}

\begin{proof}[Proof idea]
  A short
  differential geometry argument shows that $C_V$ is the
  same as $C_{V^\dual}$ in an open neighborhood 
after permuting the factors.  The equivalence of the entire bundles follows
from irreducibility.
It follows that
  $V^{\dual\dual} = V$.  See \cite[Section 2.1.5]{FP01:RuledVarieties}
  for more.
\end{proof}

\begin{proposition}\label{prop:contact-flat}
  For a homogeneous irreducible algebraic set~$V$ and a smooth point $\phi \in
  V^{\dual}$, the contact locus $V_\phi$ is flat.
\end{proposition}

\begin{proof}[Proof idea]
  We wish to show that $\pi_1(\pi_2^{-1}(H))$
  is flat for smooth $\phi \in V^\dual$.  Let us instead consider
  $\pi_2(\pi_1^{-1}(x))$ for a smooth $x \in V$.  By definition, this
  consists of all functionals vanishing on $T_x(V)$, which is a linear
  space.  By the symmetry property of Proposition~\ref{prop:double-dual},
 we can apply this
  argument in the other direction to see
  that $\pi_1(\pi_2^{-1}(\phi))$ is also a linear space.
  See \cite[Section 2.1.6]{FP01:RuledVarieties}
  for more.
\end{proof}

Proposition~\ref{prop:contact-flat} is sometimes called Bertini's
Theorem; however, there are several theorems called Bertini's Theorem,
some of which are quite different.

Note that the assumption that $\phi$
is smooth is crucial, and at singular points, the
contact locus may have different structure. For example, consider the
standard embedding in  $\RR^3$ of the 2-torus (the surface of a donut
sitting on a table).
(In this case, since we are dealing with a non-homogeneous set, a
tangent space is affine instead of linear.) At most points of the dual
variety the contact locus is the single point on the torus with that
specific affine tangent hyperplane.  A single point is flat (in the
affine sense).  But for two exceptional
hyperplanes (including the surface of the table) the contact locus is
a circle, which is not flat. These two hyperplanes
are non-smooth points of the dual variety.

\medskip

We now turn to relating the above general construction to our setting
of the measurement set~$M$.
Let $\oM$ be the Zariski closure of~$M$.  
$\oM$ is homogeneous since $M$  is closed under
multiplication by positive scalars.
Let $C_\oM$, $\oMdual$, and $\oM_\omega$ be the conormal, dual, and
contact locus constructions applied
to $\oM$.
(Here we think of $M$ as a
subset of $\RR^e$ and the equilibrium stress~$\omega$ as an element of
$(\RR^e)^*$).  We will translate freely between the stress~$\omega$
and the corresponding stress matrix~$\Omega$.

\begin{lemma}
  \label{lem:tangent-stress}
  Let $\rho \in C^d(\Gamma)$ be a framework so that $\ell(\rho)$ is
  smooth in~$M$.  Then if
  $\omega\in(\RR^e)^*$ is tangent to $M$ at
  $\ell(\rho)$, $\omega$ is an equilibrium stress for $\rho$.
  If furthermore $d\ell_\rho$ has maximal rank, then the converse
  holds: any equilibrium stress for~$\rho$ is tangent to~$M$ at $\ell(\rho)$.
\end{lemma}

\begin{proof}
Direct calculations show
\citeleft\citen{AR79:RigidityGraphsII}, p.\ 183,
\citen{Connelly05:GenericGlobalRigidity}, Lemma 2.5,
inter alia\citeright\
that the space of equilibrium 
stress vectors
at a point $\rho$,
when thought of as a subspace of $\RR^e$,
constitute the annihilator of the span of~$d\ell_{\rho}$:
\begin{equation}\label{eq:stress-annihilator}
  S(\rho) \isom \ann(\spanop(d\ell_\rho)).
\end{equation}
Finally, $\spanop(d\ell_\rho) \subset T_{\ell(\rho)}(M)$, with
equality iff $d\ell_\rho$ has maximal rank. (It
can happen that $\ell(\rho)$ is smooth even if $d\ell_\rho$ does not
have maximal rank.)  Thus $\ann(\spanop(d\ell_\rho)) \superset
\ann(T_{\ell(\rho)}(M))$, again with equality if $d\ell_\rho$ has maximal
rank.  The result follows by definition of tangency.
\end{proof}

In particular, both directions of Lemma~\ref{lem:tangent-stress}
apply at any generic framework~$\rho$.

Motivated by Lemma~\ref{lem:tangent-stress}, for $\omega \in \RR^e$
let $B^\circ(\omega)$
be the ``open contact locus'': the set of smooth points $x \in M$ so
that $\omega$ is tangent to
$M$ at~$x$.

\begin{lemma}
  \label{lem:contact-B}
  For $\rho$ a generic configuration and any $\omega \in S(\rho)$, the
  Euclidean closure of $B^\circ(\omega)$ is $B(\Omega)$.
\end{lemma}

\begin{proof}
  By Lemma~\ref{lem:tangent-stress}, if $\ell(\sigma) \in B^\circ(\omega)$,
  $\sigma \in A(\Omega)$, so
  $B^\circ(\omega) \subset B(\Omega)$.  Since $B(\Omega)$ is closed,
  we have one inclusion.

  For the other direction,
let   $A^g(\Omega)$ be the
  points in $A(\Omega)$ which are generic in~$C^d(\Gamma)$.
The set
$A^g(\Omega)$ contains $\rho$ by hypothesis and $A(\Omega)$ is a linear
space and hence smooth. 
Thus 
by  Lemma~\ref{lem:generic-dense}, $A^g(\Omega)$ is dense in $A(\Omega)$.
Thus 
$\ell(A^g(\Omega))$ is dense in 
$\ell(A(\Omega))$ which is equal to 
$B(\Omega)$.
Clearly
  $B^\circ(\omega) \superset \ell(A^g(\Omega))$.
  Thus the Euclidean closure of $B^\circ(\omega)$ contains
  the Euclidean closure of $\ell(A^g(\Omega))$, which we just argued was
  $B(\Omega)$.
\end{proof}

And now we are in a position to prove the flatness of our range space
$B(\Omega)$.

\begin{proposition}
  \label{prop:Bflat}
  For $\rho$ a generic configuration and $\omega \in S(\rho)$ so that
  $\omega$ is generic in $\oMdual$, the space $B(\Omega)$ is  flat
  (i.e., $L(\Omega)$ is the same dimension as $B(\Omega)$).
\end{proposition}

\begin{proof}
  From Lemma~\ref{lem:M-irred} $M$ is irreducible.
By Proposition~\ref{prop:contact-flat}, the contact locus
  $\oM_\omega$ is a linear space which we now identify with
  $L(\Omega)$.  
  From Lemma~\ref{lem:contact-B}, $B^\circ(\omega)$ is dense in $B(\Omega)$
  and thus $\dim B^\circ(\omega) = \dim B(\Omega) \le \dim L(\Omega)$.

  Again, due to  Lemma~\ref{lem:contact-B}, 
 $L(\Omega)$ is also the smallest linear space containing
  $B^\circ(\omega)$.
  Since, by definition,  $B^\circ(\omega)$ is contained in 
 the linear space $\oM_\omega$,
  we see that $L(\Omega) \subset \oM_\omega $
  and  $\dim L(\Omega) \le \dim \oM_\omega$.  

  Now consider an open neighborhood~$U$ of
  $\ell(\rho)$ in~$\oM$ that consists of smooth points that lie
  in~$M$ (such a neighborhood must exist since $\ell(\rho)$ is generic).  
  Then by Lemma~\ref{lem:conormal-props}, $\oM_\omega \cap U$
  consists of points $x \in M$ where $\omega$ is tangent to~$x$, i.e.,
  $\oM_\omega \cap U \subset B^\circ(\omega)$.  But $\oM_\omega \cap
  U$ is an open subset of the linear space $\oM_\omega$, 
  so $\dim \oM_\omega 
= 
\dim 
(\oM_\omega \cap U)
\le \dim
  B^\circ(\omega)$.

  Since $\dim B(\Omega) \le \dim
  L(\Omega)$, $\dim L(\Omega) \le \dim \oM_\omega$, and $\dim \oM_\omega \le
  \dim B^\circ(\Omega)$ all the inequalities must be
  equalities, and in particular $\dim B(\Omega) = \dim L(\Omega)$ as
  desired.
\end{proof}

Finally, 
we establish that for any generic framework  $\rho$, we 
can find an $\omega \in S(\rho)$ which is generic in $\oMdual$,
so that we can apply
Proposition~\ref{prop:Bflat}.

\begin{lemma}
  \label{lem:generic-omega}
  If $\rho\in C^d(\Gamma)$ is generic, then there is an $\omega \in
  S(\rho)$ such that $(\ell(\rho),\omega)$ is generic in $C_\oM$
and $\omega$  is generic inside of~$\oMdual$.
\end{lemma}

\begin{proof}
  Consider the conormal bundle $C_\oM$ with its two projections $\pi_1$
  and $\pi_2$ to $\oM$ and~$\oMdual$.  Since $\ell(\rho)$ 
is generic in~$\oM$, it is smooth in~$\oM$.
Thus, since $C_{\oM}$ is smooth above smooth points of~$\oM$,
there is a neighborhood $N$ of
$\ell(\rho)$ so that $\pi_1^{-1}(N)$ is a smooth semi-algebraic set.
Applying Lemma~\ref{lem:preimage-generic} to the restriction of 
$\pi_1$ to $\pi_1^{-1}(N)$
guarantees a point
  $(\ell(\rho),\omega)\in \pi_1^{-1}(\ell(\rho))$ that is generic in 
$\pi_1^{-1}(N)$ and thus also in $C_\oM$.
By
  Lemma~\ref{lem:image-generic}, $\omega\in\oMdual$ is
  also generic.  This $\omega$ is an equilibrium stress for $\rho$ by
  Lemmas~\ref{lem:conormal-props} and~\ref{lem:tangent-stress}.
\end{proof}

\subsection{Identifying the contact locus}
\label{sec:ident-contact}

We now digress briefly to give a somewhat more
explicit description of the space $L(\Omega)$.  This is not necessary
for any of our proofs, but may aid in understanding.

\begin{definition}
  For a set $S \subset \EE^n$ and $d \in \NN$, 
let 
$\chord^d(S)$ be the $d$'th
chord set of $S$, the union of all simplices with
$d$ vertices, all in $S$.  For instance, when $d=2$ we add
chords connecting all pairs of points in~$S$.
(For a variety~$V$, the Zariski closure of $\chord^2(V)$ is 
the secant variety $\sec(V)$, where we add complete lines instead of
segments.)
\end{definition}
\begin{lemma}
  \label{lem:B-chord}
  The space $B(\Omega)$ is $\chord^d(\ell(K(\Omega)))$.
\end{lemma}

\begin{proof}
  The squared edge lengths of a framework~$\rho$ are computed by
  summing the squared edge
  lengths of each coordinate projection of~$\rho$.  Since $A(\Omega)$
  is $K(\Omega)^d$, $B(\Omega)$ is
  the $d$-fold Minkowski sum of $\ell(K(\Omega))$ with itself.
  Because $\ell(K(\Omega))$ is
  invariant under scaling by positive reals, this Minkowski sum
  coincides with $\chord^d(\ell(K(\Omega)))$.
\end{proof}

\begin{definition}\label{def:dot-prod-D}
For $a, b \in C^1(\Gamma)$,
define a ``dot product'' $\langle a, b\rangle \in \RR^e$ by
\begin{equation}\langle a, b\rangle(\{w,u\}) \coloneqq 
(a(w)-a(u))\cdot (b(w)-b(u))\end{equation} 
for each edge $\{u,w\} \in \Edges$.
Define the dot product space $D(\Omega)$ to be the linear span of
$\langle a, b\rangle$ for all $a,b \in K(\Omega)$.
\end{definition}

\begin{lemma}
  \label{lem:L-is-D}
  The space $L(\Omega)$ is the same as $D(\Omega)$.
\end{lemma}

\begin{proof}
The image
$\ell(K(\Omega))$ is contained in $D(\Omega)$, as it is just the space of
all $\langle a, a\rangle$ for $a\in K(\Omega)$.  
Moreover, the linear
span of $\ell(K(\Omega))$ is in fact equal to $D(\Omega)$, as $\langle a,
b\rangle = \frac{1}{2}(\langle a+b, a+b\rangle - \langle a, a\rangle -
\langle b, b\rangle)$.
From Lemma~\ref{lem:B-chord}, the linear span of 
$\ell(K(\Omega))$ 
is the same
as the linear span of
$B(\Omega)$,
which is $L(\Omega)$ by definition.
\end{proof}

A priori, if $K(\Omega)$ is $k$-dimensional, the dimension of
$D(\Omega)$ could be as big as
$\binom{k}{2}$.  (There are $\binom{k+1}{2}$ 
dot products between 
vectors forming a basis of $K(\Omega)$.  
However, dot products with the vector of all
ones vanish because of the subtraction in the definition of
$\langle a, b\rangle$.)  On the other hand, the dimension of $B(\Omega)$ is only
$kd - \binom{d+1}{2}$.
When $\Gamma$ has a minimal stress kernel realized by $\Omega$,
$k=d+1$ and $\dim(B(\Omega))$ agrees with
this estimate for $\dim(D(\Omega))$.
The crux of Proposition~\ref{prop:Bflat} is that in fact, even
when $\Gamma$ does not have minimal stress kernel, if $\Omega$ is generic
in ~$\oMdual$ then
$B(\Omega)$ is
a flat space, and thus $D(\Omega)$ is in fact  only of dimension 
$kd - \binom{d+1}{2}$. This means that
there must be some linear dependence between the dot products
defining $D(\Omega)$.

\subsection{The map}
\label{sec:map}
We now turn to our main map.

\begin{definition}
  Let~$\Omega$ be 
  an equilibrium stress matrix of some framework in $C^d(\Gamma)$.
  The map $\ol$ is the restriction of $\ell$ to a map between
  the spaces $A(\Omega)/\Eucl(d)$ and~$L(\Omega)$.
\end{definition}

In this section,
we show that,
for any generic framework $\rho$ and stress $\Omega \in S(\rho)$ so
  that $\Omega$ is generic in $\oMdual$,
the map $\ol$ has a well-defined mod-two degree which is~$0$, and
then deduce the theorem.
(Recall that the existence of such a generic $\Omega$ is guaranteed 
by Lemma~\ref{lem:generic-omega}.)

First we 
recall a few standard elementary properties of regular values, 
proper maps, and degrees.  Because our maps are not always locally injective or surjective, we 
use a slightly generalized definition of regular value
and  version of Sard's  theorem.

\begin{definition}
  Let $f : X \to Y$ be  a smooth map between smooth manifolds $X$
  and~$Y$. Let $r$ be the maximal rank
of the linearization $df_x$ of~$f$ for any $x \in X$. 
We say $x \in X$ is a \emph{regular point} if 
$df_x$ has rank $r$.
We say $y \in Y$ is a
\emph{regular value} if 
every $x \in f^{-1}(y)$ is a regular point.
Otherwise, $y$ is called a \emph{critical value}.
\end{definition}
\vspace{0pt} % Some spacing bug
\begin{theorem}[Sard~\cite{sard}]~\label{thm:sard}
Let $f : X \to Y$ be  a smooth map between smooth manifolds, and let
$r$ be the maximal rank
of $df_x$ for  $x \in X$.
The critical values of $f$ have $r$-dimensional measure 0.
\end{theorem}

\begin{proposition}\label{prop:alg-sard}
Let $X$ and $Y$ be semi-algebraic sets that are manifolds and $f$ be a 
polynomial map with maximal rank~$r$, all defined over the field~$\kk$.
Then the critical values form a semi-algebraic subset of~$Y$, defined
over $\kk$, of dimension less than~$r$.
In particular, all generic points of
$f(X)$ are regular values.
\end{proposition}

\begin{proof}
  See, e.g., \cite[Theorem 9.6.2]{bcr}.
\end{proof}

\begin{definition}
  \label{def:proper} A \emph{proper map} $f : X \to Y$ between
  topological spaces is a continuous map so that the inverse image of
  a compact set is compact.
\end{definition}

Examples of proper maps include the identity and any map from a
compact space.  
For our purposes, we will need the following:

\begin{lemma}
\label{lem:lisproper}
The length measurement map
associated to any connected graph,
once we quotient the domain by the group $\Trans(d)$ of translations,
is proper,
as is its restriction to $A(\Omega)$ for any equilibrium stress matrix~$\Omega$.
That is, the maps
\begin{align*}
  \ell &: C^d(\Gamma)/\Trans(d) \to \RR^e\\
  \ell &: A(\Omega)/\Trans(d) \to \RR^e
\end{align*}
are both proper.
\end{lemma}
\begin{proof}
A
compact
subset $P \subset\RR^e$ is bounded, so gives a bound on the edge lengths.
This in turn gives a bound on how far any vertex in $\Gamma$ can be
from some fixed base vertex, so $\ell^{-1}(P)$ is bounded in
$C^d(\Gamma)/\Trans(d)$ or the subspace $A(\Omega)/\Trans(d)$.  Since
$\ell^{-1}(P)$ is also closed, it is compact.
\end{proof}

There is a notion of \emph{degree} of proper maps between manifolds of
the same dimension.  The following theorem is standard, and is
typically proved using homology \cite[Theorem 8.12, inter alia]{Spivak79:DifferentialGeometry}.

\begin{theorem}
  \label{thm:degree}
  If $X$ and $Y$ are manifolds of the same dimension, with $Y$
  connected, and $f : X \to Y$ is a proper map, then there is an
  element $\deg f$ in $\ZZ/2$, invariant under proper isotopies
  of~$f$.  If $X$, $Y$, and $f$ are all
  smooth, then the degree is equal to the number of preimages of any
  regular value, taken modulo 2.
\end{theorem}

We will want to compute the degree in a case when~$X$ is not quite a manifold,
but rather has singularities of codimension~$2$.  We use the following
version.

\begin{corollary}
  \label{cor:degree-sing}
  If $X$ is a smooth stratified space with
  singularities of codimension at least~$2$, $Y$ is a smooth,
  connected manifold of the same dimension as~$X$, and $f : X \to Y$
  is a proper,
  smooth map, then there
  is an element $\deg f$ in $\ZZ/2$, invariant under proper isotopies
  of~$f$.  The degree is equal to the
  number of preimages of any regular value in~$Y$, taken modulo~$2$.
\end{corollary}

Here by a ``regular value'' we mean a point in~$Y$ so that every
preimage is a smooth, regular point in~$X$.  The condition that the
singularities have codimension at least~$2$ is crucial; otherwise, for
instance, the inclusion of the interval $[0,1]$ in $\RR$ would qualify, and
the degree is obviously not invariant there.

\begin{proof}
  Let $X^{\textrm{sing}}$ be the set of singular points of~$X$.  Let
  $Y' = Y\setminus f(X^{\textrm{sing}})$, and let $X' = f^{-1}(Y')$.
  Then $X'$ and $Y'$ are both smooth manifolds by construction.  By
  Lemma~\ref{lem:proper-excision} below, the restriction of $f$ to a
  function from $X'$ to $Y'$ is still proper.  
The smooth image of the stratified space~$X^{\textrm{sing}}$ is itself
a stratified set,
with dimension no bigger: $\dim(f(X^{\textrm{sing}})) <
\dim(X^{\textrm{sing}})$.
Therefore
  $Y'$ is still connected, as we have removed a subset of
  codimension at least~$2$ from~$Y$.  We can therefore apply
  Theorem~\ref{thm:degree} to
  find $\deg f$ as the degree of the restriction from $X'$ to $Y'$.
\end{proof}

\begin{lemma}[Excision]
  \label{lem:proper-excision}
  If $f:X \to Y$ is a proper map and $Y' \subset Y$ is an arbitrary subset,
  let $X' = f^{-1}(Y')$.  Then the restriction of $f$ to $X'$, $f' :
  X' \to Y'$, is proper.
\end{lemma}

\begin{proof}
  For any compact subset $P$ of $Y'$, $P$ is also compact as a subset
  of $Y$.  Since
  $f'^{-1}(P)$ is the same as $f^{-1}(P)$, it is compact.
\end{proof}

More generally, there is an integer-valued degree for maps between
oriented manifolds.  Our domain space is not in general oriented, so
we only have a mod-two degree, but that is enough for us.

Applying these results to our map $\ol$, we see that:

\begin{lemma}
\label{lem:even}
Suppose that $\Gamma$, 
a graph with $d+2$ or more vertices,
is generically locally rigid in $\EE^d$ and
does not have minimal stress kernel in~$\EE^d$, 
$\rho$ is a generic framework in~$C^d(\Gamma)$,
and $\omega \in S(\rho)$ is generic in $\oMdual$.
Then the resulting map
$\ol$ from $A(\Omega)/\Eucl(d)$ to $L(\Omega)$ has a mod-two degree
of~$0$.
\end{lemma}
\begin{proof}
By Lemma~\ref{lem:lisproper} the map $\ol$ is proper.  (Local rigidity
implies that $\Gamma$ is connected.)

Since $\Gamma$ does not have a minimal stress kernel and $\rho$ is generic, 
$\dim(K(\Omega)) > d+1$ for all  $\Omega \in S(\rho)$.
Thus by Proposition~\ref{prop:sing-codim-Arho}, 
the domain is a
smooth stratified space with singularities of codimension $2$ or
greater.  
By Proposition~\ref{prop:Bflat}, $L(\Omega)$ (the range of $\ol$) has
the same dimension as
$B(\Omega)$, which is the same dimension as the domain by local rigidity.
So by 
Corollary~\ref{cor:degree-sing}
it has a well defined degree mod-two.

Since all squared edge lengths in the image of $\ell$ are positive,
any point in $L(\Omega)$ with some negative edge lengths has no
preimage in $A(\Omega)/\Eucl(d)$ and is automatically regular. Hence
the mod-two degree must
be~$0$.
\end{proof}

For Lemma~\ref{lem:even} to be useful, we must see that $\ell(\rho)$
is a regular value of~$\ol$.
\begin{lemma}
\label{lem:regval}
Suppose that $\Gamma$, 
a graph 
with $d+2$ or more vertices,
is generically locally rigid in $\EE^d$, and 
$\rho$ is a generic framework in $C^d(\Gamma)$.
Then
$\ell(\rho)$ is a regular value of both $\ell$ and~$\ol$.
\end{lemma}

\begin{proof}
Since $\ell(\rho)$ is a generic point of $M$, by Sard's Theorem
$\ell(\rho)$ is a regular value of~$\ell$, proving the first part of
the statement.  At each preimage~$\sigma$ of a
regular value, $d\ell_\sigma$ is of maximal rank, and so $\sigma$ is
infinitesimally rigid. In particular, $\sigma$ must have a
$d$-dimensional affine span, and hence cannot have a 
non-trivial stabilizer in $\Eucl(d)$.
Then by Proposition~\ref{prop:sing-codim-Arho} every point in
$\ol{}^{-1}(\ell(\rho))$ is
a smooth point in $A(\Omega)/\Eucl(d)$.

Since $d\ell_\sigma$ is of maximal rank, it is
injective from the
tangent of $C^d(\Gamma)/\Eucl(d)$ to $\RR^e$.  The map $d\ol_\sigma$ is the
restriction of $d\ell_\sigma$ from the tangent space to $C^d(\Gamma)/\Eucl(d)$
to a subspace.
Restriction preserves injectivity of the linearization of a smooth
map, so $\ell(\rho)$ is also a regular value for~$\ol$.
\end{proof}

The regularity of $\ell(\rho)$ with respect to $\ol$ 
can also be proved by
using Lemma~\ref{lem:generic-in-A}.

And now we are in position to complete the proof of our main Theorem.

\begin{proof}[Proof of Theorem~\ref{thm:deficient-flexible}]
By Lemma~\ref{lem:generic-omega},
for any generic $\rho$ there is
an $\omega \in S(\rho)$ that is generic in $\oMdual$. 
Choose the associated equilibrium stress matrix $\Omega$
to define $\ol$.
From Lemma~\ref{lem:even}, the mod-two degree of $\ol$ is $0$.
From Lemma~\ref{lem:regval}, $\ell(\rho)$ is a regular value of $\ol$.
Thus there must be an even number of points in 
$\ol^{-1}(\ell(\rho))$.
Since there is at least one point in the preimage,
namely $[\rho]$ itself, there must be another point in
$A(\Omega)/\Eucl(d)$, and thus an incongruent framework in
$C^d(\Gamma)$ that has
the same edge lengths as~$\rho$. 
\end{proof}

%%% Local Variables: 
%%% mode: latex
%%% TeX-master: "main"
%%% End: 

\section{Examples}
\label{sec:examples}
Here we
give a few examples of how the various spaces constructed in
Section~\ref{sec:proof} work out in practice.  We consider both cases
where Theorem~\ref{thm:deficient-flexible} applies to show that
generic frameworks are not globally rigid, and cases where it does
not apply.

Throughout this section, $\rho$ is a generic framework for whichever
graph we are looking at and $\omega \in S(\rho)$ is chosen to be
generic in $\oMdual$.  (In particular, $\dim(K(\Omega)) = k_{\min}$.)

\subsection{Locally flexible graphs}
\label{sec:locally-flex-examp}

\begin{figure}
  \[
  \mathgraph{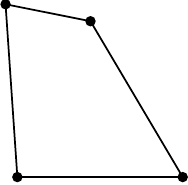}\qquad\qquad\qquad\mathgraph{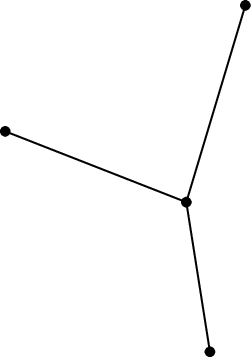}
  \]
  \caption{Generic locally flexible frameworks in~$\EE^2$}
  \label{fig:locally-flex}
\end{figure}
The first case we consider is that of graphs that are generically
locally flexible in $\EE^d$, as in Figure~\ref{fig:locally-flex}.

The space $A(\Omega)$ contains all affine transforms of any flex
of~$\rho$, and possibly more. This can be seen by using
the main
argument from~\cite{Connelly05:GenericGlobalRigidity}: any flex~$\sigma$
of $\rho$ maps to $\ell(\rho)$ in $M$, which is smooth since $\rho$
is generic.  We therefore have $\spanop(d\ell_\sigma) \subset
T_{\ell(\rho)} M = \spanop(d\ell_\rho)$.  In particular, any
equilibrium stress for~$\rho$ is also an equilibrium stress for $\sigma$, so
$\sigma \in A(\Omega)$.

A simple example of such a graph is the tripod graph in $\EE^2$
as on the right in Figure~\ref{fig:locally-flex}.
In this case, the only 
equilibrium 
stress vector is the zero stress, and $A(\Omega)$ contains
all embeddings of $\Gamma$. Note that in this case, $L(\Omega)$ is 
three dimensional, which is the same size as one would find for a 
generically globally rigid graph.

\subsection{Not redundantly rigid}
\label{sec:not-redund-rigid-examp}

\begin{figure}
  \[
  \mathgraph{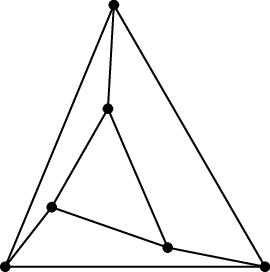}\qquad\qquad\mathgraph{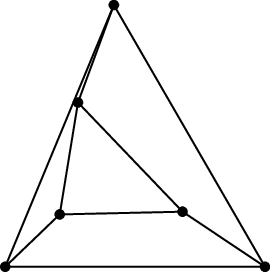}\qquad\qquad
  \mathgraph{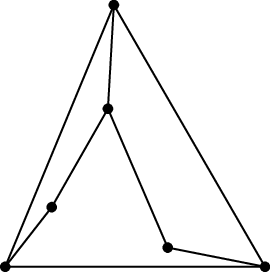}
  \]
  \caption{A framework that is not redundantly rigid in~$\EE^2$, with
    graph given by a triangular prism.  From left to right, we have a
    framework~$\rho$, the global flex guaranteed by the theorem, and
    the locally flexible framework~$\rho'$ obtained by deleting a
    suitable edge.}
  \label{fig:not-redund-rigid}
\end{figure}
Here we suppose that $\rho$ is generically locally rigid, but becomes
not generically locally rigid when
some edge $\{i,j\}$  is removed from~$\Gamma$, as shown for example in
Figure~\ref{fig:not-redund-rigid}.
In other words, we suppose that the graph can
be flexed, leaving all edge lengths unchanged except for edge $\{i,j\}$,
which is changed.
These graphs are  not
generically globally 
rigid by a result of
Hendrickson~\cite{Hendrickson92:ConditionsUniqueGraph}; we will see
how this works in our setting.
For a generic framework~$\rho$, we can find, in the span of
$d\ell_{\rho}$, a vector
with all zeros except for the entry associated with
edge $\{i,j\}$ (since by definition there is a flex changing this edge
length and no others).
Since the equilibrium 
stress vectors annihilate $\spanop(d\ell_{\rho})$, we conclude that 
any equilibrium stress $\Omega$ has a zero at entry $\{i,j\}$.
(In fact, a locally rigid generic framework is redundantly rigid at
edge $\{i,j\}$ iff
it has an equilibrium stress with a non-zero entry at that edge.)

Define $\Gamma'$ to be the graph obtained by deleting this edge, and
$\rho'$ to be the resulting framework.  
Since every equilibrium stress for generic frameworks has
zeroes at entry $\{i,j\}$, the equilibrium
stresses for $\Gamma$ and $\Gamma'$ are
the same; that is, $\overline{M(\Gamma)}^\dual = \overline{M(\Gamma')}^\dual$.
In particular, we may think of our chosen~$\Omega$ as an equilibrium stress for
$\rho'$, and it is then generic in $\overline{M(\Gamma')}^\dual$.
By the analysis of generically 
locally flexible frameworks
in Section~\ref{sec:locally-flex-examp}, $A(\Omega)$ contains all
affine transforms of flexes of $\rho'$.  In terms of $\rho$, we would
say that $A(\Omega)$ contains all affine transforms of frameworks of
$\Gamma$ that have the same edge lengths as $\rho$ except on the
edge $\{i,j\}$.

It is also instructive to picture the relationship between 
the measurement set $M(\Gamma)$ in $\RR^e$ 
and the measurement set $M(\Gamma')$ in $\RR^{e-1}$.
Because the framework is not redundantly rigid,
$M(\Gamma)$ projects onto
$M(\Gamma')$ (by forgetting one coordinate in $\RR^e$), and the fiber
over a generic point $\ell(\rho')$ contains an interval.
Accordingly, $L(\Omega)\subset \RR^e$ for $\Gamma$
is the direct sum of the corresponding $L(\Omega)\subset\RR^{e-1}$ for $\Gamma'$
and the vector
in the direction of the missing edge.

\subsection{Graphs that are not $(d+1)$-connected}
\label{sec:d-conn-examp}
\begin{figure}
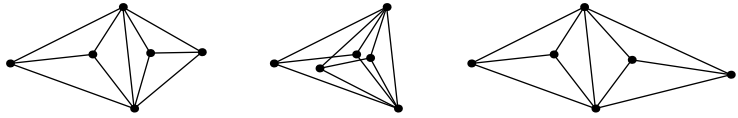

  \[
  \mathgraph{stresses-30}\qquad\mathgraph{stresses-31}
  \qquad\mathgraph{stresses-32}
  \]
  \caption{From left to right, a framework~$\rho$ that is not
    3-connected in $\EE^2$ (and
    therefore not globally rigid), its global flex, and
    another framework in $A(\Omega)$.}
  \label{fig:d-conn}
\end{figure}
A graph is not $(d+1)$-connected if there is a set of at most $d$ vertices
whose removal leaves a disconnected graph.
These graphs are  not 
generically globally rigid in
$\EE^d$~\cite{Hendrickson92:ConditionsUniqueGraph}.
If the separating set has $d$ vertices, they span an ``interface''
hyperplane whose removal
splits the graph into two halves, and we can reflect one half across
the interface.
If the separating set has fewer vertices, the
graph is not even generically 
locally rigid, as we can pivot the two halves around
the interface.

Let $\Gamma'$ be one of the two halves of the graph, together with the
vertices and edges in the separating set, and let $\rho'$ be the
corresponding framework.  Any equilibrium stress~$\Omega$ in $S(\rho)$
induces a
non-equilibrium stress on $\rho'$.  The resulting forces on the
vertices of $\Gamma'$ must have no net translational and rotational
effect (as is true for any stress, equilibrium or not).  On the other
hand, these forces are non-zero only on the $d$ vertices of the
interface.  It is then easy to see that the forces must lie entirely within
the plane of the interface.

Therefore the forces on the interface vertices are unchanged by
affinely 
squashing one half towards the interface while leaving the other half
alone.
Thus
$A(\Omega)$ includes frameworks where one applies two different affine
transforms to
the two halves, with the constraint that they agree on the interface.
In particular, this includes the global flex that entirely folds one
side across the interface.

\subsection{Bipartite graphs}
\label{sec:k_5.5-examp}

\begin{figure}
  \[
  \includegraphics[width=2.5in]{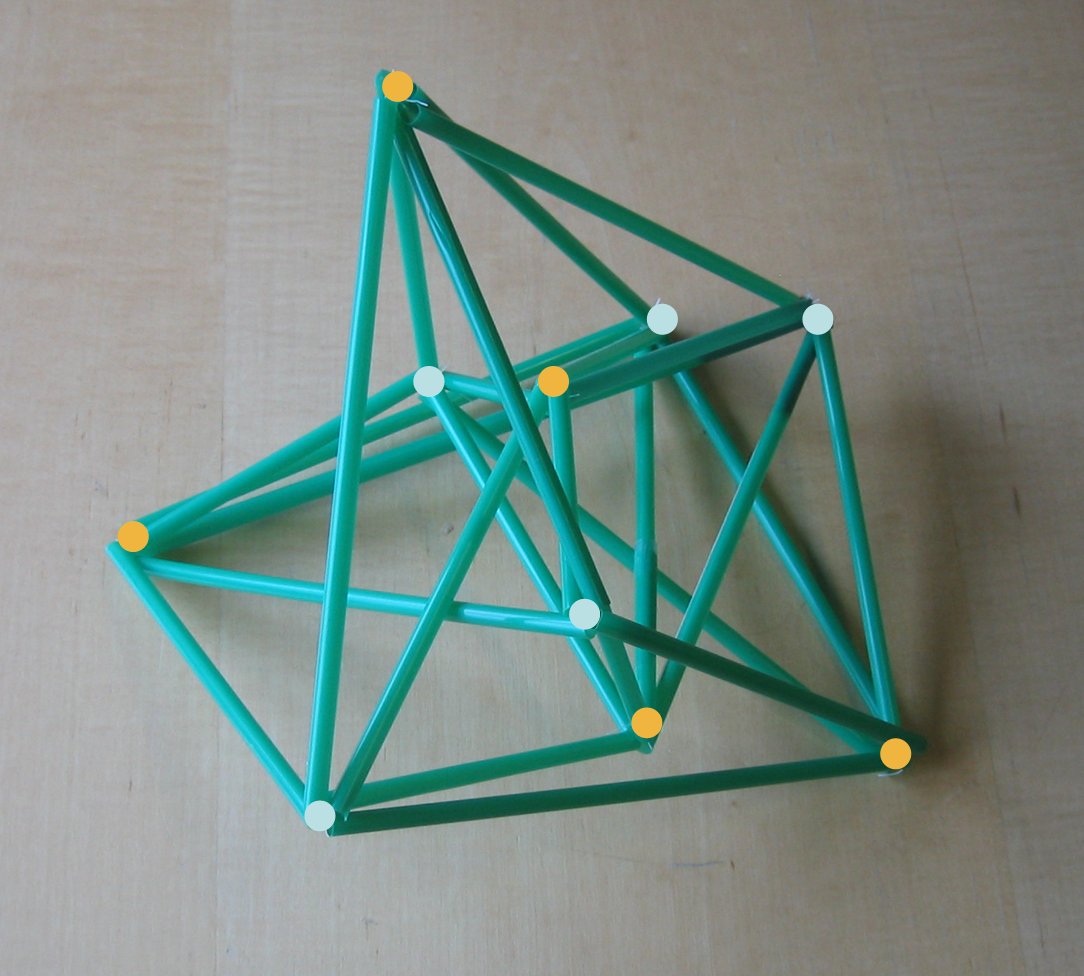}
  \]
  \caption{A framework of $K_{5,5}$ in 3 dimensions}
  \label{fig:K55}
\end{figure}

  We next consider certain complete bipartite frameworks which (as
  Connelly showed) satisfy Hendrickson's
  conditions but are still not generically globally
  rigid\cite{Connelly91:GenericGlobalRigidity}.
  The first example is $K_{5,5}$ in $\EE^3$, as in
  Figure~\ref{fig:K55}, which we focus on first.
  Let the vertices be $v_1,\dots,v_{10}$ so that
  the each of $v_1,\dots,v_5$ is connected to each of
  $v_6,\dots,v_{10}$, 25 edges in all.  We have that $C^3(K_{5,5})$ is
  30-dimensional.  Since this graph is generically locally rigid, the measurement
  set~$M$ has dimension $30 - \dim(\Eucl(3)) = 24$ inside $\RR^{25}$.
  We therefore have a 1-dimensional space $S(\rho)$ of equilibrium 
stress vectors  at a
  generic framework, which we now identify, following Bolker and
  Roth~\cite{BR80:BipartiteRigid}.
  (In fact, Bolker and
    Roth went the other way: they used the identification of the space
    of equilibrium stresses to identify
    which bipartite graphs are generically locally rigid, including
    $K_{5,5}$ in $\EE^3$.)
  We will write $v_1, \dots,
  v_{10}$ for the 10 points in $\EE^3$ at the generic framework~$\rho$.

  The 5 points $v_1,\dots,v_5$ in $\EE^3$
  satisfy an affine relation
  \begin{equation}
    \begin{split}
      a_1 v_1 + \cdots + a_5 v_5 &= 0\\
      a_1 + \cdots + a_5 &= 0
    \end{split}
  \end{equation}
  where not all $a_i$ are~$0$.  Since the
  $v_i$ are generic, this
  relation is unique up to scale.
  Similarly there is a unique up to scale affine
  relation between $v_6,\dots,v_{10}$:
  \begin{equation}
    \begin{split}
      b_6 v_6 + \cdots + b_{10} v_{10} &= 0\\
      b_6 + \cdots + b_{10} &= 0.
    \end{split}
  \end{equation}
  Now form vectors $\vec a =
  (a_1,\dots,a_5)$ and $\vec b = (b_6,\dots,b_{10})$ and consider the
  block matrix
  \begin{equation}
  \Omega =
  \begin{pmatrix}
    0 & \vec a^T \vec b\\
    \vec b^T \vec a & 0
  \end{pmatrix}.
  \end{equation}
  (That is, the $(i,j)$ entry of $\Omega$ is $a_i b_j$ or $a_j b_i$ if
  either of these is defined, and $0$ otherwise.)  Then it follows that
  $\Omega$ is an  equilibrium 
 stress matrix and so spans the space $S(\rho)$.

 In particular, since $\Omega$
 visibly has rank~$2$, $\dim K(\Omega) = 8$.  The basis vectors are
 \begin{enumerate}
 \item a vector of 3 coordinate projections of $v_1$ through $v_5$ and
   five $0$'s in the last 5 positions,
 \item a vector of five $1$'s and five $0$'s, and
 \item five $0$'s in the first 5 positions followed by the 3
   coordinate projections of $v_6$ through $v_{10}$,
 \item a vector of five $0$'s and five $1$'s.
 \end{enumerate}

  The space $A(\Omega)$ is also easy to identify: it is the
  $24$-dimensional space of
  configurations $w_1,\dots,w_{10}$, where $w_1,\dots,w_5$ are an
  affine transform of $v_1,\dots,v_5$ and $w_6,\dots,w_{10}$ are an
  affine transform of $v_6,\dots,v_{10}$, but the two affine
  transforms need not be the same.  This guarantees that
  $w_1,\dots,w_5$ satisfy the affine relation given by $\vec a$ and
  $w_6,\dots,w_{10}$ satisfy the affine relation given by $\vec b$,
  and so these configurations satisfy the  equilibrium 
stress matrix~$\Omega$.

It follows that $L(\Omega)$ is an $18$-dimensional space.  We can think
of it as the space $D(\Omega)$ as in Definition~\ref{def:dot-prod-D}.

More generally, Bolker and Roth showed that a generic framework of
$K_{n,m}$ in $\EE^d$ with $n,m \ge d+2$ and $n+m = \binom{d+2}{2}$ is
locally rigid and has equilibrium stresses only of the form
\begin{equation}
  \Omega =
  \begin{pmatrix}
    0 & X\\
    X^T & 0
  \end{pmatrix}
\end{equation}
where $X$ is a $n \times m$ matrix in which each column is an affine
linear relation among $v_1, \dots, v_n$ and each row is an affine linear
relation among $v_{n+1}, \dots, v_{n+m}$.  In particular, the
framework is redundantly rigid (since there are 
equilibrium stresses with non-zero
entries one every edge).  The rank of $X$ is at most $\min(n,m) - d -
1$, the rank of $\Omega$ is twice as large, and
\begin{equation}
  k_{\min} = n + m - \rank(\Omega) = \abs{n-m} + 2(d+1) > d+1
\end{equation}
so (as Connelly stated) these graphs are not generically globally
rigid.

\subsection{Generically globally rigid graphs}
\label{sec:glob-rigid-examp}

Finally we consider a generic framework~$\rho$ in $\EE^d$ that
\emph{does} satisfy
Connelly's condition, so Theorem~\ref{thm:deficient-flexible} does not
apply and the framework is generically globally rigid.  It is instructive 
to see
where the proof breaks.  In this case, for every $\Omega\in  S(\rho)$
that is generic in
$\oMdual$, $K(\Omega)$
is only
$(d+1)$-dimensional.  
The space $A(\Omega)$ is then just 
$\af(\rho)$, 
the space of
affine transforms of~$\rho$.

Let $\rho_i$ be the projection of $\rho$ onto the $i$'th coordinate.
Then by Lemma~\ref{lem:L-is-D}, $L(\Omega)$ is spanned by the
$\binom{d+1}{2}$ different vectors $\langle
\rho_i, \rho_j\rangle$.  These
are all linearly independent,
since,
as shown 
in~\cite[Proposition~4.3]{Connelly05:GenericGlobalRigidity},
whenever Connelly's condition is satisfied,
the edges of~$\rho$ do not lie on a
\emph{conic at infinity}.  (A conic at infinity is, in our language, a
linear
dependence among the $\langle \rho_i, \rho_j\rangle$.)

Moreover, by \cite[Proposition 4.2]{Connelly05:GenericGlobalRigidity},
whenever
the edges of~$\rho$ do not lie on a
conic at infinity, 
$\ell$ is injective on $\af(\rho)/\Eucl(d)$.
Thus in our case,
$B(\Omega)$ is homeomorphic to $A(\Omega)/\Eucl(d)$.
In particular, $B(\Omega)$ has a boundary (as it is a proper subset of
$L(\Omega)$) and so must
$A(\Omega)/\Eucl(d)$.
This boundary may be thought of as a singularity of codimension~$1$ in
the domain, and thus
our degree argument does not apply.

More explicitly,
$\ell(\af(\rho))$
is equivalent to a
\emph{semidefinite cone}, as we now explain.
Let
$\ls(\ell(\af(\rho)))$ 
be the linear span of the image space $\ell(\af(\rho))$.
Use the basis $\langle \rho_i, \rho_j\rangle$ to define an
isomorphism $\phi$ from $\ls(\ell(\af(\rho)))$ 
to the space of symmetric 
$d\times d$ matrices, sending $\langle\rho_i, \rho_j\rangle$ to
$e_{ij}$, where
\[
e_{ij} =
\begin{cases}
  \text{matrix with a $1$ at position $(i,i)$
and $0$'s elsewhere} & i = j\\
  \text{matrix with a
$1/2$ at positions $(i,j)$ and $(j,i)$ and $0$'s elsewhere} & i \ne j.
\end{cases}
\]

\begin{lemma}
  For a framework~$\rho$ whose edges do not lie on a conic at infinity
and $\phi$ defined as
  above, 
$\phi(\ell(\af(\rho)))$ 
is the set of positive semidefinite matrices.
\end{lemma}

\begin{proof}
  Since $\ell$ is invariant under translation of~$\rho$, it suffices
  to consider $\phi(\ell(U\rho))$, where $U =
  (u^i_j)_{i,j=1}^d$ is a linear map (rather than a general affine
  map).  We then find
  \begin{equation}
    \begin{split}
    \phi(\ell(U(\rho))) &= \phi\left(\sum\nolimits_j 
      \left\langle \sum\nolimits_i u^i_j \rho_i,
        \sum\nolimits_i u^i_j \rho_i \right\rangle\right)\\
%    &= \phi\left(\sum\nolimits_{j,i,i'} u^i_j u^{i'}_j
%      \langle\rho_i,\rho_{i'}\rangle\right)\\
    &= \sum\nolimits_{j,i,i'} u^i_j u^{i'}_j e_{i i'}\\
    &= U^T U.
    \end{split}
  \end{equation}
  The $d\times d$ positive semidefinite matrices are exactly those
  matrices that can be written as $U^T U$ for a $d \times d$
  matrix~$U$.
\end{proof}

In particular, when the graph has minimal stress kernel and $\rho$ and
$\Omega$ are
generic, $A(\Omega)= \af(\rho)$, and thus
$B(\Omega)$ and $A(\Omega)/\Eucl(d)$ are isomorphic to the semidefinite cone
in
$\phi(L(\Omega))$.
The number of preimages in
$A(\Omega)/\Eucl(d)$ of a point in $L(\Omega)$ is $1$
for points inside the positive semidefinite cone.  (At non-generic
points, there may be more preimages in $C^d(\Gamma)/\Eucl(d)$ that
are not in $A(\Omega)$).
The number of preimages is
$0$ for points
outside this cone (such as points in $\RR^e$ with any negative
coordinate).  Thus the degree is definitely not independent of the
generic target point we consider.

%%% Local Variables: 
%%% mode: latex
%%% TeX-master: "main"
%%% End: 

\section{Shared stress kernels and the rank of the Gauss map}
\label{sec:shared-stress}

In this section we give an alternate
characterization of global rigidity, in terms of the \emph{shared
  stress kernel}.

\begin{definition}
  The \emph{shared stress kernel} of a framework~$\rho$ is the
  intersection of all stress
  kernels: $K(\rho) \coloneqq
  \bigcap_{\Omega
    \in S(\rho)} K(\Omega)$.  
  It is isomorphic to the space of frameworks of $\Gamma$ in $\EE^1$ which
  satisfy \emph{all} the  equilibrium stress matrices in $S(\rho)$ (and maybe
  more).  Let the \emph{shared stress nullity} $k_\shared(\Gamma, d)$ be
  $\dim(K(\rho))$ for a generic framework~$\rho$.
\end{definition}

\begin{remark}
  As for $k_{\min}$, the shared stress nullity $k_{\shared}$ is
  independent of the generic framework~$\rho$.  One way to see this is
  to note that the intersection in the definition of $k_\shared$ is
  achieved by considering a sufficiently large finite number of 
  equilibrium stress matrices
  $\Omega_i \in S(\rho)$.  The intersection of these kernels is also the
  kernel of a larger matrix ${\Omega'} =
  (\Omega_1\,\Omega_2\cdots\Omega_n)$
  obtained by concatenating the~$\Omega_i$.  But ${\Omega'}$ 
  can be constructed by algebraic functions (as in the proof of 
  Lemma~\ref{lem:generic-stress-rank})
  so its rank
  is maximal at generic points by Lemma~\ref{lem:poly-mat}.

  Alternately, we can use the equivalence with the Gauss rank in
  Proposition~\ref{prop:mgm} below to show the same fact.
\end{remark}

\begin{remark}
  By definition, $k_{\min} \ge k_{\shared}$.  It frequently happens
  that $k_{\min} = k_\shared$, but it is not true in general that
  the analogous numbers are the same.  For instance, the matrices
  \[
  \begin{pmatrix}
    0 & a & b \\
    a & 0 & 0 \\
    b & 0 & 0
  \end{pmatrix}
  \]
  for $a,b \in \RR$ all have rank 2 (and so a 1-dimensional kernel);
  but the intersection of their kernels is zero.  A similar example
  occurs in practice for bipartite graphs; for frameworks of $K_{n,m}$ in
  $\EE^d$ with $n,m \ge d+2$, $n+m = \binom{d+2}{2}$, and $n \ne m$,
  it follows from the description of the equilibrium stresses recalled in
  Section~\ref{sec:k_5.5-examp} that
  \begin{align}
    k_{\min} &= \abs{n-m} + 2(d+1)\\
    k_{\shared} &= 2(d+1) < k_{\min}.
  \end{align}
  The second equation follows from the fact that the only frameworks that
satisfy all of the stresses of $\rho$ must be related to $\rho$  by
an affine transform on each of the two sides of the bipartite graph.
  (For instance, $K_{7,8}$ in $\EE^4$ satisfies these constraints, and
  has $k_{\min} = 11$ and $k_{\shared} = 10$.)
\end{remark}

Although $k_{\min}$ and $k_{\shared}$ are different in general, either
one can be used in a test for global rigidity, so if $k_{\min} =
v-d-1$, then $k_{\shared} = v-d-1$ as well:

\begin{theorem}\label{thm:shared-rigid}
  A graph~$\Gamma$ with $d+2$ or more vertices is generically globally
  rigid in $\EE^d$ if and only if $k_{\shared}(\Gamma, d) = d+1$.
\end{theorem}

One direction of Theorem~\ref{thm:shared-rigid} is a corollary of
Theorem~\ref{thm:deficient-flexible}, as $k_{\shared} \le k_{\min}$.
For the other direction, we strengthen Connelly's argument, using
the same proof.

\subsection{The Gauss map and its rank}
\label{sec:gauss-map}

Before continuing we point out that there is 
an alternate way of interpreting
the shared stress kernel,
by using the rank of maps.

\begin{definition}
  The \emph{rank} of an algebraic or rational map between irreducible algebraic
  varieties is
  rank of its linearization at generic points in the domain. (This is
  also the maximum rank of its linearization.)  Equivalently, it is
  the dimension of the image semi-algebraic set.
\end{definition}
Thus
Theorem~\ref{thm:locally-rigid} says that a 
graph is generically locally rigid in $\EE^d$
if
and only if 
the rank of the map~$\ell$ is
``the largest it can be''. 

Now assume that $\Gamma$ is generically locally rigid in~$\EE^d$, and
let $t = vd - \binom{d+1}{2}$ be the dimension of the measurement
set~$M$.
Consider now the \emph{Gauss map}, the map that takes each smooth 
point of $M$ to its tangent space.  We consider
the tangent space as
a point in the space of $t$-dimensional planes in $\RR^e$, the Grassmannian
$\Gr(t,e)$, so we get a map $G: M \dashrightarrow \Gr(t,e)$.  (The
dashed arrow indicates that the map is a rational map: it is not
defined at non-smooth points,
which must lie on some subset of dimension less than $t$; formally,
$G$ is defined on a Zariski open subset of~$M$.)
We may then consider the 
rank of~$G$; in fact, this rank is closely related to~$k_\shared$.

\begin{proposition}
  \label{prop:mgm}
The maps $G$ and $G \circ \ell$ both have rank $vd - k_\shared d$.
\end{proposition}

Before proving this proposition, we make some preliminary definitions,
paralleling the earlier definitions for the minimal stress kernel.

\begin{definition}
  For $\rho \in C^d(\Gamma)$, define
  \begin{align}
    K(\rho) &= \bigcap_{\Omega \in S(\rho)} \ker(\Omega)\\
    A(\rho) &= \{\,\sigma \in C^d(\Gamma) \mid S(\rho) \subset S(\sigma)\,\}.
  \end{align}
\end{definition}

Note that if $\Omega \in S(\rho)$, by definition $K(\rho) \subset
K(\Omega)$ and $A(\rho) \subset A(\Omega)$.
We can think of $K(\rho)$ as the 1-dimensional frameworks that satisfy
every stress in $S(\rho)$.  

\begin{lemma}
  \label{lem:A-K-rel-shared}
    The space $A(\rho)$ is isomorphic to $K(\rho)^d$.
\end{lemma}
\begin{proof}
  By Lemma~\ref{lem:A-K-rel}, we have
  \begin{equation}
  A(\rho) = \bigcap_{\Omega \in S(\rho)} A(\Omega) =
  \bigcap_{\Omega \in S(\rho)} K(\Omega)^d = K(\rho)^d.\qedhere
  \end{equation}
\end{proof}

\begin{proof}[Proof of Proposition~\ref{prop:mgm}]
  Fix a generic framework $\rho\in C^d(\Gamma)$.  Observe that by
  Lemma~\ref{lem:tangent-stress}, the equilibrium stresses for~$\rho$
  are the hyperplanes tangent to~$\ell(\rho)$.  By
  Lemma~\ref{lem:tangent-stress} again, these hyperplanes are also
equilibrium
  stresses for any point~$\sigma$ so that $\ell(\sigma)$ is smooth and
  such that $\ell(\sigma)$ is in the fiber of $G$ through
  $\ell(\rho)$.  (But note that $\rank(d\ell_\sigma)$ need not be
  maximal.)  Thus the fiber of $G \circ \ell$ through $\rho$ is
  contained in $A(\rho)$.

  Now define the space
  $A^{\circ}(\rho)$ to be the set of frameworks $\sigma \in A(\rho)$
  so that $\rank(d\ell_\sigma)$ is maximal and $\ell(\sigma)$ is a
  smooth point of~$M$.  Both conditions are algebraic, and $\rho$
  satisfies both of them, so $A^{\circ}(\rho)$ is a non-empty Zariski
  open subset of $A(\rho)$.  In particular,
$\overline{A^{\circ}(\rho)} = A(\rho)$.

Meanwhile 
from Lemma~\ref{lem:tangent-stress} all points in 
$A^{\circ}(\rho)$
must be in the fiber of $G \circ \ell$.

Thus the fiber has dimension $\dim
A^{\circ}(\rho) =
\dim A(\rho)$.
We therefore have
\begin{equation}
\dim A(\rho) =
\dim C^d(\Gamma) - \rank (G \circ \ell).
\end{equation}
Since $G$ is defined on an open subset of the image of $\ell$, the
image of $G$ is the same dimension as the image of $G \circ \ell$ and
$\rank (G \circ \ell) = \rank G$.
Putting these together with Lemma~\ref{lem:A-K-rel-shared} we conclude
\begin{equation}
\rank (G \circ \ell) = \rank G  = \dim C^d(\Gamma) - \dim A(\rho)
  = vd - k_\shared d.
\qedhere
\end{equation}
\end{proof}

In light of Proposition~\ref{prop:mgm}, since $k_\shared \geq d+1$, we have
$\rank G \leq vd- (d+1)d$.
Thus we can introduce the following terminology.

\begin{definition}
  \label{def:mgm}
  A graph~$\Gamma$ has \emph{maximal Gauss rank in $\EE^d$}
if 
 $\rank G  = vd- (d+1)d$.
\end{definition}

In particular,
a graph $\Gamma$
has a minimal stress kernel
iff it has maximal Gauss rank, and by
Theorems~\ref{thm:connelly-rigid}, \ref{thm:deficient-flexible}, and
\ref{thm:shared-rigid} both are
equivalent to $\Gamma$ being generically globally rigid in~$\EE^d$.

\subsection{Strengthening Connelly's proof}
\label{sec:strength-connelly}

We now prove the missing direction of Theorem~\ref{thm:shared-rigid},
mildly strengthening Connelly's proof of
Theorem~\ref{thm:connelly-rigid}, closely following his argument.

\begin{proof}[Proof of Theorem~\ref{thm:shared-rigid}]
  Let $\Gamma$ be a graph with $d+2$ or more vertices with
  $k_\shared(\Gamma, d) = d+1$.  Let $\rho$ be a generic framework in
  $C^d(\Gamma)$, and let $\sigma$ be another framework with the same
  edge lengths.  By \cite[Proposition
  3.3]{Connelly05:GenericGlobalRigidity}, there are neighborhoods
  $U_\rho$ and $U_\sigma$ of $\rho$ and $\sigma$ and a diffeomorphism
  $g: U_\sigma \to U_\rho$ with $\ell(g(x)) = \ell(x)$ for $x \in
  U_\sigma$.  In particular,
  \begin{equation}
  \spanop(d\ell_\sigma) = \spanop(d\ell_\rho).
  \end{equation}
  It follows that $S(\sigma) = S(\rho)$ (since both are equal to
  $\ann(\spanop(d\ell_\rho))$).  Thus $\sigma$ is in $A(\rho)$, which
  is $K(\rho)^d$ by Lemma~\ref{lem:A-K-rel-shared}.  By hypothesis,
  $K(\rho)$ consists of the coordinate
  projections of $\rho$ and the vector $\vec 1$ of all $1$'s, so
  $\sigma$ is an affine image of~$\rho$.  (This is the analogue of
  \cite[Theorem 4.1]{Connelly05:GenericGlobalRigidity}.)

  As in the proof of \cite[Theorem
  1.3]{Connelly05:GenericGlobalRigidity}, for each vertex~$v_i$ we can
  find a stress in $S(\rho)$ with some non-zero stress-value on an edge
  connected to the vertex.  (Otherwise we could freely move~$v_i$
  without changing the space of equilibrium
stresses, and get a stress kernel that
  is too large.  Here we use the fact that there are at least $d+2$
  vertices.)  Therefore each vertex has degree at least $d+1$, as the
  vectors $v_i - v_j$ for $\{i,j\} \in \Edges$ are linearly dependent.
  Then by \cite[Proposition
  4.3]{Connelly05:GenericGlobalRigidity}, the edges of $\rho$ do not
  lie on a conic at infinity, so by \cite[Proposition
  4.2]{Connelly05:GenericGlobalRigidity} any affine image of $\rho$
  with the same edge lengths as $\rho$ is actually congruent
  to~$\rho$.  Therefore $\sigma$ is congruent to~$\rho$, as
  desired.
\end{proof}

\subsection{Secant varieties and Zak's theorem}
\label{sec:secant-zak}

We add a few remarks about the algebraic geometry setting of these
results.  Both generic local rigidity and global rigidity relate to
the geometry of the map $\ell$ and the measurement
set.  In both cases the appropriate rigidity property is
characterized in terms of purely local behavior of the measurement
set.  For the generic local rigidity case, this is not surprising,
but for the generic global rigidity
case, it is a little unusual to be able to characterize a global property (the
number of pre-images of a certain map) simply in terms of local behavior
(the rank of the Gauss map).

Let us consider what happens as we vary the dimension~$d$
while keeping the graph fixed.
For the case $d=1$, the map $\ell$ is
essentially an arbitrary quadratic map and the measurement set is
of dimension $v-1$ if the graph is connected.
Call this measurement set $M_1$, and denote by
$\overline{M_1}$ its
Zariski closure, which we call the measurement variety.
 (For the remainder of this section, we are
interested in varying the ambient dimension while keeping the graph
fixed, so we add a subscript.)
In $d$ dimensions, the squared edge lengths are the sum of the squared
edge lengths in each dimension separately.
Thus the measurement variety~$\overline{M_d}$ in~$\EE^d$
is the closure of the $d$-fold Minkowski sum of $\overline{M_1}$.  Since
$\overline{M_1}$ is a homogeneous variety, the Minkowski sum coincides with the
secant variety, and so
\begin{equation}
\overline{M_d} = \overline{\ell(C^d(\Gamma))} = \overline{\sec^d(\overline{M_1})}
\end{equation}
where $\sec^d(\overline{M_1})$ is the $d$-fold secant
variety of $\overline{M_1}$ (i.e., the join of $d$ copies of $\overline{M_1}$).

For general homogeneous varieties, the expected dimension of the join of two
varieties is the sum of the two dimensions and so the expected dimension
of the $d$-fold secant variety of a variety~$V$ is $d$ times the
dimension of~$V$.  If the dimension of the secant variety is less than the
expected dimension the secant variety is said to be \emph{degenerate}.
In our case, because $\overline{M_1}$ is the closure of the image of a
quadratic map there is
automatically some degeneracy and its secant variety $\overline{M_2}$ has
dimension at least one less than expected: $\dim \overline{M_2} \le
2\dim\overline{M_1} - 1$.  If our variety $\overline{M_1}$ were
smooth or almost\hyp smooth,
Zak's Theorem on
Superadditivity~\cite{Zak85:LinearSystHyperplanes,Fantechi90:SupearditivitySecant}
would imply that this degeneracy
propagates, and given this first degeneracy, the $d$-fold secant $\overline{M_d}$
would have dimension at most $d\cdot \dim \overline{M_1} - \binom{d}{2}$.
Although the conditions of Zak's Theorem do not hold for a general
graph~$\Gamma$, the
resulting bound is still
exactly the dimension of $\overline{M_d}$ given by
Theorem~\ref{thm:locally-rigid} for graphs that are
generically locally rigid in $\EE^d$. Thus generically locally rigid
graphs are
those where the secant variety is \emph{minimally degenerate} in an
appropriate sense.

Generic global rigidity instead involves another notion of degeneracy:
the \emph{dual defect} of a variety.  For a
general homogeneous variety $V$ in $\RR^n$, the dual variety $V^*$ is
$(n-1)$-dimensional.  (In this case the generic contact locus~$V_\phi$
is a single line through the origin.)  The dual defect is defined to
be $n-1 - \dim V^*$:
the difference between the expected dimension of $V^*$ and the actual
dimension.  Again, in our setting $\overline{M_d}$ automatically has a
dual defect of $\binom{d+1}{2}-1$ since it is the secant variety of the
image of a quadratic map.
Theorems~\ref{thm:connelly-rigid} and~\ref{thm:deficient-flexible} can
be interpreted as saying
that the graph is generically globally rigid if and only if the dual defect
is this minimum.

Alternately, we can work with the \emph{Gauss defect}.  For a general
homogeneous variety~$V\subset \RR^n$, the rank of the
Gauss map~$G: V \dashrightarrow \Gr(t,n)$ is~$\dim V-1$.  (In fact, this
happens whenever the projectivization $\PP V$ is smooth.)
Define the Gauss defect of~$V$ to be $\dim V - 1 - \rank G$.  In general,
the Gauss defect is less than or equal to the dual defect
\cite[Section 2.3.4]{FP01:RuledVarieties}.  Again, because
$\overline{M_d}$ is the secant variety of the image of a quadratic
map, the Gauss defect is at least $\binom{d+1}{2}-1$.  By
Theorem~\ref{thm:shared-rigid} the graph is generically globally rigid
if and only if the Gauss defect is this minimum.  (In particular, in
this case the Gauss defect equals the dual defect.)

There is much literature on varieties with degenerate secant
varieties or Gauss maps,
although most of the work focuses on cases which are
maximally, rather than minimally, degenerate
\cite{LV84:TopicsGeomProj, FL81:ConnectivityAppl,
  Zak93:TangentsSecants}.
From this point of view, Connelly's
example~\cite{Connelly91:GenericGlobalRigidity}
of $K_{5,5}$ (which is not
generically globally rigid in 3 dimensions) is interesting. In this
case the variety $\overline{M_1}$ is the cone on the Segre
embedding of $\RR\PP^4 \times \RR\PP^4$ in $\RR\PP^9$, and Connelly's
result is about the degeneracy of the Gauss map of the third secant variety
of this Segre embedding.

%%% Local Variables: 
%%% mode: latex
%%% TeX-master: "main"
%%% End: 

\section{Complexity of the algorithm}
\label{sec:complexity}

Implicit in Theorems~\ref{thm:connelly-rigid}
and~\ref{thm:deficient-flexible} are a deterministic and
a randomized algorithm.
In this section, after briefly describing the 
deterministic algorithm (which is not very efficient),  we will prove
that testing whether a graph is generically globally
rigid in $\EE^d$ is in RP: that is, there is a polynomial-time randomized
algorithm that will answer
``no'' all the time if the graph is not generically globally rigid in $\EE^d$,
and will answer ``yes''
at least half the time if the graph is generically globally rigid in $\EE^d$.
Hendrickson
sketched an argument that testing for generic \emph{local} rigidity is in RNC
(which is contained in RP)  \cite[Section
2.2.2]{Hendrickson92:ConditionsUniqueGraph}, and later sketched an 
randomized algorithm for determining whether a graph has a minimal stress
kernel in $\EE^d$ \cite[Section 2.3.1]{Hendrickson95:MoleculeProblem}.
Here we will provide some more details on
these algorithms, why they work, and what the necessary bounds on the
size of the inputs are.

\subsection{Deterministic algorithm}
\label{sec:determ-alg}

First we sketch the deterministic algorithm.  Let
$t\coloneqq vd-\binom{d+1}{2}$ be the rank of the rigidity matrix for
infinitesimally rigid frameworks, and let $s \coloneqq v - d - 1$ be the
maximal rank of an equilibrium
stress matrix for generic
globally rigid frameworks.

\begin{algorithm}[Deterministic check, global]
\label{alg:deterministic-global}
To check if a graph~$\Gamma$ with at least $d+2$ vertices
is generically globally rigid in~$\EE^d$,
take a $d$-dimensional framework with distinct coordinates, treated as
independent formal symbols, and
compute its rigidity matrix, again symbolic.  By Gaussian elimination
over the field of fractions of the coordinates, compute the
rank of this matrix and kernel of its transpose.
(The resulting kernel will be a ratio of polynomials in the symbolic inputs,
and the answers will be valid
for a numerical framework
as long as none of the expressions in the denominators are zero.)
If the rank of the rigidity matrix is not
$t$, output the answer ``no'', as the graph is not even
generically locally rigid.  Otherwise
take a generic equilibrium 
stress vector, i.e., a symbolic linear combination of a
basis of the kernel of the transposed rigidity matrix, convert it to
an equilibrium stress
matrix, and compute its rank (again symbolically).  If this rank is
$s$, output the answer ``yes'', otherwise output the
answer ``no''. 
\end{algorithm}

While this algorithm is guaranteed to give the correct answer, it is
likely to be very slow, as the number of terms in the polynomials we need to
manipulate can grow rapidly.

\subsection{Efficient randomized algorithms}
\label{sec:rand-alg}

For a more efficient version, we turn to
a similar randomized algorithm.  

\begin{algorithm}[Randomized check, local]
  \label{alg:random-local}
  To check if a graph~$\Gamma$ with at least $d+1$ vertices
  is generically locally rigid in $\EE^d$,
  pick a framework~$\rho$ in $C^d(\Gamma)$ with integer coordinates
  randomly chosen from $[1,N]$ for some
suitably large~$N$ (to be made precise below).
Compute the rank of the rigidity
  matrix representing 
$d\ell_\rho$.  If this rank is less than $t$, output ``no'',
  otherwise output ``yes''.
\end{algorithm}

\begin{algorithm}[Randomized check, global]
\label{alg:random-global}
To check if a graph~$\Gamma$ with at least $d+1$ vertices
is generically globally rigid in~$\EE^d$, proceed as follows.
First, if the number of edges, $e$, is less than $t$ 
output ``no'' (as the graph cannot even be generically locally
rigid with so few edges), otherwise continue.

Next
pick a
framework~$\rho$ in $C^d(\Gamma)$
with integer coordinates randomly chosen from $[1,N]$ 
for some
suitable~$N$.

The next step is to pick
one equilibrium 
stress vector~$\omega$ from $S(\rho)$ in a suitably random way.
Recall that an equilibrium 
 stress vector is a vector
in the annihilator of $\spanop(d\ell_{\rho})$, i.e., in the kernel
of $(d\ell_\rho)^T$.
We find $\omega$ by setting up 
and solving an appropriate
linear system $E \omega = b$ that extends the condition that $\omega \in
\ker((d\ell_\rho)^T)$.
To create the linear system,
create a matrix~$H$ of $e-t$ random row vectors 
in $\RR^e$, with coordinates integers chosen from $[1,N]$.
Append these random rows
to the transpose of the rigidity matrix 
to obtain an $e+\binom{d+1}{2}$ by~$e$ matrix~$E(\rho, H)$.
Compute the rank of~$E(\rho, H)$. If this rank is less than $e$, output ``no'',
for either the rank of the rigidity matrix is less than $t$
(for instance if $\Gamma$ is not generically locally rigid),
or 
$H$ contains some linear
dependence with the rest of~$E$,
otherwise continue.
(Note that, in constructing an RP algorithm, it is safe to
output ``no'' in unfavorable cases, as long it is not done too often.)

Define $b$ to be a vector in $\RR^{e+\binom{d+1}{2}}$ 
which is all zeros except for a single
entry of $1$ corresponding to one of the rows in~$H$. (If $e=t$, 
there will be no such entry,  $b$ will be all zeros and the only $\omega$
will be the zero vector.) 
Now solve the linear system $E \omega = b$, which must have exactly one
solution, denoted $\omega(\rho,H)$.

Finally, convert $\omega(\rho,H)$ into an  equilibrium 
 stress matrix $\Omega(\rho,H)$, and compute its
rank.  If the rank is $s$, output the answer
``yes'', otherwise output the answer ``no''.
\end{algorithm}

\begin{remark}
  In the randomized algorithm as described above, the linear algebra
  computations are all done with explicit integer matrices.  This can
  be done exactly in polynomial time using $p$-adic
  techniques~\cite{Dixon82:ExactSolLinear}.  (That is, first solve the
  equations modulo a prime~$p$, then lift to a solution modulo~$p^2$,
  and so forth.)  However, since the algorithm is already randomized,
  the added complexity of this exact algorithm is unnecessary, and the
  computations can be done more simply by reducing modulo a suitably
  large prime~$p$.  See Proposition~\ref{prop:primes-est} below for a
  concrete estimate of how large the primes need to be.
\end{remark}

To analyze the probability of the algorithm giving
a false negative output, which depends on how large the parameter~$N$ is,
our main tool is the \emph{Schwartz-Zippel Lemma}
\cite{DL78:ProbRmkTesting,Zippel79:ProbAlgSparse,Schwartz80:FastProbAlg}.

\begin{lemma}[Schwartz, Zippel, DeMillo, Lipton]
  \label{lem:schwartz-zippel}
  Let $\kk$ be a field and $P \in \kk[x_1,\dots,x_n]$ be a non-zero
  polynomial of degree~$d$ in the~$x_i$.  Select $r_1,\ldots,r_n$ uniformly at
  random from a finite subset~$X$ of~$\kk$.  Then the probability than
  $P(r_1,\dots,r_n)$ is~$0$ is less than $d/\abs{X}$.
\end{lemma}

We also use the following basic principle:

\begin{lemma}\label{lem:poly-mat}
Let $M(\pi)$ be a matrix whose entries are polynomial functions 
with rational coefficients in the
variables $\pi \in   \RR^n$.
Let $r$ be a rank achieved by some $M(\pi_0)$.
Then $\rank (M(\pi))\geq r$ 
for all points  $\pi$ that are generic in $\RR^n$.
More precisely, if the entries of $M(\pi)$ are polynomials  of degree
$g$, 
then there is a polynomial $P$ on $\RR^n$ of degree $g\cdot r$ with rational
coefficients 
so that $P(\pi)=0$ for any $\pi$
where the rank of $M(\pi)$ is less than $r$.
\end{lemma}

\begin{proof}
The rank of the $M(\pi)$  is
less than~$r$ iff the determinants of all of the
$r\times r$ submatrices
vanish.  
Let $\pi_0\in\RR^n$ be a choice of parameters so $M(\pi_0)$ has
rank~$r$.  Then there is an $r \times r$ submatrix $T(\pi_0)$ of
$M(\pi_0)$ with non-zero determinant.  Thus $\det(T(\pi))$ is
a non-zero polynomial of~$\pi$ of degree $g\cdot r$.
For any $\pi$ with $\rank(M(\pi)) < r$, this
determinant 
must vanish. 
Thus, any such $\pi$ cannot be generic.
We can therefore take $P(\pi) = \det(T(\pi))$.
\end{proof}

We now apply these concepts to  the simpler case of
local rigidity for illustration.

\begin{proposition}\label{prop:random-local-RP}
  Algorithm~\ref{alg:random-local} never returns a false ``yes''
  answer, and returns a false ``no'' answer with probability bounded
  above by $t/N$.
  In particular, checking for generic local rigidity is in RP.
\end{proposition}

\begin{proof}

If the graph is not generically locally
rigid, then from Lemma~\ref{lem:poly-mat}
the algorithm 
must  always observe  a rank which is less than~$t$. Thus it will
never give a false ``yes'' answer.

Suppose instead that $\Gamma$ is generically locally rigid in~$\EE^d$.
The entries of the rigidity
matrix are linear in~$\rho$,
so from Lemma~\ref{lem:poly-mat}
any $\rho$ which returns a false ``no'' answer must be at a zero
of a polynomial of degree $t$ in our random variables.
But by the Schwartz-Zippel Lemma, $\rho$ is a zero 
 with probability less than $t/N$, as desired.
If we choose $N > 2t$, we will give a false ``no'' answer less than
$1/2$ the time.

The algorithms for finding
ranks, etc., take time polynomial in $\log N$ (as well as the size of the
matrix), so run in overall polynomial time.  
The one-sided nature of the errors
means that the problem is in RP (rather than, say, BPP).
\end{proof}

We now return to the global case.
Recall that for a non-singular $n \times n$ matrix~$M$,
\begin{equation}\adj(M) = \det(M)M^{-1},\end{equation}
where $\adj M$ is the \emph{adjugate
 matrix} of $M$, the conjugate of the cofactor matrix of $M$.
In particular,
$\adj{M}$ is a polynomial in~$M$ of degree $n-1$.
This essentially means that, up to a global scale,
the entries of the inverse of a matrix~$M$ are polynomials in the
entries of~$M$.  We will apply this to the matrix~$E(\rho,H)$ from
Algorithm~\ref{alg:random-global} to see that the equilibrium 
 stress matrix $\Omega(\rho,H)$ found by the
algorithm
can be though of as polynomials in the entries of~$E(\rho,H)$.
As the entries of $E(\rho,H)$ are polynomial functions in the coordinates of
the framework $\rho$ and the random rows~$H$, and matrix rank is
invariant to scaling, this will put us in position to apply
Lemma~\ref{lem:poly-mat}.

There are a few technical issues that must be dealt with
in this analysis. 
Most importantly, the matrix $E(\rho,H)$ has more rows than columns. We will
deal with this by looking at an appropriate 
square submatrix of $E(\rho,H)$, as follows.

\begin{lemma}\label{lem:generic-stress-rank}
  Let $\Gamma$ be a graph, and
  let $t'$ be the maximal value of $\rank d\ell_\rho$ over all
  frameworks~$\rho$ of~$\Gamma$ in $\EE^d$.
  (For generically
  locally rigid graphs $t' = t$.)  Let $U$ be the set of $\rho$ where this
maximal rank is achieved. 
Let $s'$ be the maximal value of 
$\rank(\Omega)$ over all 
stress matrices of all $\rho$ in $U$.
Then 
for all 
generic $\rho$, we have
$\max_{\Omega \in S(\rho)} \rank(\Omega) =s'$.

More precisely, there is a polynomial $P$ of degree $e+s'(e-1)$ in the
variables $(\rho, H)$, where $H$ is a set of additional rows as
in Algorithm~\ref{alg:random-global}, so that if $(\rho, H)$ is not a
zero of~$P$ then the stress matrix $\Omega(\rho,H)$ is well-defined
and has rank $s'$.
\end{lemma}

Here, as in
Algorithm~\ref{alg:random-global}, for each $\rho$
we add a set~$H$ of $e-t'$ additional rows
  in $\RR^e$ to the transpose of the rigidity matrix of $\rho$ to define a
  matrix~$E(\rho,H)$, and find an
equilibrium 
stress vector~$\omega(\rho, H)$ by
solving the  linear system $E\omega = b$ 
where $b \in \RR^e$ is a vector of all zeroes except for a
single~$1$ in one of the positions of a row in~$H$ (if any).  This is
then converted to a stress matrix $\Omega(\rho, H)$.  The assertion
that $\Omega(\rho, H)$ is well-defined means that this linear system
has a unique solution.  (In particular, this implies that $\rho \in U$.)

\begin{proof}[Proof of Lemma~\ref{lem:generic-stress-rank}]
Let $\rho_0\in U$ be a framework with an equilibrium stress $\omega_0$
so that the corresponding stress matrix $\Omega_0$ has rank $s'$.
Find a set~$H_0$ of additional rows so that $E(\rho_0,H_0)$ has
rank~$e$ and $E(\rho_0,H_0)\omega_0 = b$.  Let $\hE(\rho,H)$ be an $e
\times e$ submatrix of $E(\rho,H)$ so that $\hE(\rho_0,H_0)$ is
invertible.  ($\hE$ necessarily consists of $t'$ rows from the
transpose of the rigidity matrix and all the rows of $H$.)  Define
$\hb$ similarly, let $\tomega(\rho,H) \coloneqq \adj(\hE)\hb$, and let
$\tOmega(\rho,H)$ be
the associated stress matrix.

By Lemma~\ref{lem:poly-mat}, the rank of $\hE(\rho,H)$ is equal to its
maximum value~$e$ at all points $(\rho,H)$ that are not zeros of a
polynomial $P_1(\rho, H) = \det \hE(\rho,H)$, which has degree~$e$.
Moreover, when $P_1(\rho, H) \ne 0$, the linear equation defining
$\Omega$ has a unique solution and the adjugate
matrix $\tOmega(\rho, H)$ is a scalar multiple of $\Omega(\rho,H)$.  
In
particular 
we have assumed $(\rho_0, H_0)$
is not a zero of $P_1$ and thus 
$\tOmega(\rho_0, H_0)$ has rank~$s'$.  By
Lemma~\ref{lem:poly-mat} again, the rank of $\tOmega(\rho,H)$ 
is less than $s'$ only
at the zeros of a non-zero polynomial $P_2(\rho, H)$ of degree
$s'(e-1)$ (as the entries of $\tOmega(\rho,H)$ have degree $e-1$ in
$(\rho,H)$).

For any generic $\rho$, 
there must be some generic point
$(\rho,H)$. 
At such a generic 
$(\rho,H)$ 
$\tOmega(\rho,H)$ and $\Omega(\rho,H)$ 
must have rank not less than~$s'$. 
Meanwhile we have supposed that rank of such $\Omega(\rho,H)$ 
is never larger than
$s'$, and thus have proven our claim.

In particular, let $P = P_1 \cdot P_2$.  Then if $(\rho, H)$ is not a
zero of $P$, the linear equation $E\omega = b$ has a unique solution
(as $P_1 \ne 0$) and the resulting stress matrix $\Omega(\rho,H)$ has
rank~$s'$ (as $P_2 \ne 0$).
\end{proof}

\begin{remark}
  In fact, since the matrix rank drops only at the zeros of an 
algebraic equation, 
we can prove the stronger statement that all generic
  $\Omega \in S(\rho)$ of all generic $\rho$ have $\rank(\Omega)=s'$.
  Here we consider $S(\rho)$ as defined over the smallest extension of~$\QQ$
  containing~$\rho$.
\end{remark}

The following theorem is a more precise version of
Theorem~\ref{thm:random-global-RP-1}.

\begin{theorem}\label{thm:random-global-RP}
  Algorithm~\ref{alg:random-global}
  never returns a false ``yes''
  answer, and returns a false ``no'' answer with probability bounded
  above by $ve/N$.
  In particular, checking for generic global rigidity in~$\EE^d$ is in RP.
\end{theorem}

\begin{proof}
First we suppose that $\Gamma$  is not generically
globally rigid in~$\EE^d$.
If the rigidity matrix does not have rank $t$, then
the algorithm correctly outputs ``no''. If the rigidity matrix does
have rank $t$, by
Lemma~\ref{lem:generic-stress-rank}, the maximal observed rank cannot 
be 
larger than the rank of any  equilibrium 
stress matrix of any generic framework,
which is less than~$s$ by supposition.
Thus the algorithm will correctly output ``no''.

Next we suppose
that $\Gamma$  is  generically
globally rigid in~$\EE^d$ and 
bound the probability of a false
negative.

Let $(\rho, H)$ be a set of variables that are not a zero of the
polynomial~$P$ given by Lemma~\ref{lem:generic-stress-rank}.  Such
points cannot give a false negative: the rank of~$E(\rho,H)$ is its maximum
value, $e$, so we do not give a false negative in
constructing~$\omega$, and $\rank(\Omega)=s$, so again we do not
output a false
negative. 

So if $\Gamma$ is
generically globally rigid in $\EE^d$, 
we are guaranteed to  obtain the correct answer if we pick 
$\rho$ and~$H$
that miss
the 
zeros of $P$, a polynomial of degree $e+s(e-1)$,
which is less than $ve$.  
The Schwartz-Zippel Lemma, Lemma~\ref{lem:schwartz-zippel}, then gives the
bounds stated on the error
probability.
Moreover,
if we choose $N > 2ve$, we will 
give a false ``no'' answer less than half the time.  

The running
time for such an $N$ is polynomial in the input size, so testing for
generic global
rigidity is in RP.
\end{proof}

Algorithm~\ref{alg:random-global} and the proof of
Theorem~\ref{thm:random-global-RP} are structured as they
are because Lemma~\ref{lem:generic-stress-rank} is \emph{false}
without the condition that $\rank d\ell_\rho$ take its maximal value.
For instance, if a framework $\rho$ lies in a lower-dimensional subspace
of~$\EE^d$, the rank of a generic~$\Omega \in S(\rho)$ may be
greater than $v-d-1$, so the universal upper bound motivating
Definition~\ref{def:deficient} may not hold.

\subsection{Estimating the primes}
\label{sec:estimating-primes}

For the reader's convenience, we now estimate some bounds on the size
of the primes necessary for the algorithms to work with high
probability.  These bounds can very likely be improved substantially.

\begin{proposition}
  \label{prop:primes-est}
  Suppose that the computations in Algorithm~\ref{alg:random-local}
  or~\ref{alg:random-global} are done modulo a prime~$p$ 
  then neither
  algorithm will return a false ``yes'' answer.
  If the prime is chosen
  chosen
  uniformly from a
  set of primes~$\Primes$ with each prime larger than $N$, 
  and if $\abs{\Primes} > 4t$ and
  $N > 4t$, then Algorithm~\ref{alg:random-local} produces a false
  ``no'' answer with probability bounded above by $1/2$.
  Similarly, if $\abs{\Primes} > 4ve$ and
  $N > 4ve$, then Algorithm~\ref{alg:random-global} produces a false
  ``no'' answer with probability bounded above by $1/2$.
\end{proposition}

\begin{remark}
  By the Prime Number Theorem, this proposition tells us we need to
  consider primes up to about $8t\ln(4t)$ (for local rigidity) or
  $8ve\ln(4ve)$ (for global rigidity) in order to get a sufficiently
  large set~$R$.
\end{remark}

\begin{proof}
  Doing the computations modulo a prime can only make the computed ranks
  drop, so 
  as before, we cannot return a false ``yes'' answer.

  If the graph
  is generically locally (or globally) rigid, there is an integer
  polynomial $P_\loc$ (or $P_\glob$) so that if our random framework
  $\rho$ (or pair $(\rho, H)$) is not a zero of $P_\loc$ (or
  $P_\glob$) we do not return a false ``no'' answer.
  The Schwartz-Zippel Lemma, Lemma~\ref{lem:schwartz-zippel}, works
  without change modulo~$p$.  The
  only differences are that
  \begin{itemize}
  \item the prime~$p$ must be large enough so that there are enough
    distinct values modulo~$p$ for the lemma to be useful, and
  \item the polynomial~$P$ not be zero modulo~$p$.
  \end{itemize}
  To take care of the first point, we require that the prime be
  larger than the chosen value~$N$ in the Schwartz-Zippel Lemma.

  To take care of the second point, we will give an upper bound~$B$
  for the sum of the absolute values of the coefficients of~$P$.  If
  we pick a prime~$p$ larger than $B$ (or in fact larger than the
  absolute value of any single coefficient of $P$) then $P$ will be
  guaranteed to be non-zero modulo~$p$.  However, we can in fact
  randomly pick a smaller prime from a suitable collection, and still
  guarantee that with high probability $P$ will be non-zero
  modulo~$p$.  In particular, suppose that we have a collection of
  primes~$\Primes$ so that any product of at least $\abs{\Primes}/4$
  primes in the collection is larger than~$B$.  Then, any subset of
  size at least $\abs{\Primes}/4$ must have some prime such that $P$ is
  non-zero modulo $p$.  (Otherwise, $P$ would have to be zero modulo
  this product, but this this product is larger than $B$.)  Thus $P$
  is zero modulo $p$  for at most $1/4$ of the primes in $\Primes$.
  Then if we pick $p$ at random from~$\Primes$, the chance that $P$ is
  zero modulo~$p$ is at most $1/4$.

  We now estimate the sum of the coefficients of~$P$ in the two
  cases.  For local rigidity, $P_\loc$ is a determinant of a $t \times t$
  submatrix of the rigidity matrix.  As such, it is a sum of $t!$
  terms, with each term a product of $t$ factors, and with each factor
  having coefficient-sum
  equal to~$2$.  Using the estimate $n! < n^n$, we therefore find
  \begin{equation}
    \ln B_\loc \le \ln (2^t \cdot t!) < t(\ln t + \ln 2).
  \end{equation}
    For global rigidity, we separately estimate bounds $B_1,
  B_2$ for the two factors $P_1, P_2$ making up $P_\glob$ from the
  proof of Lemma~\ref{lem:generic-stress-rank}.  We find
  \begin{equation}
    \ln B_1 \le \ln(2^{e_\Gamma} \cdot e_\Gamma!) <
       e_\Gamma(\ln e_\Gamma + \ln 2).
  \end{equation}
  (Here we write $e_\Gamma$ for the number of edges in~$\Gamma$ to
  avoid confusion with the base of natural logarithms.)
  Note that $P_2$ is an $s\times s$ determinant of $\tOmega$, which
  itself is made of $(e_\Gamma-1)\times (e_\Gamma-1)$ determinants
  of~$E$.  We therefore have
  \begin{equation}
    \ln B_2 \le \ln (s! (2^{e_\Gamma} e_\Gamma!)^s) <
     s \ln s + s e_\Gamma \ln 2 + s e_\Gamma \ln e_\Gamma.
  \end{equation}

  Now suppose, in the local rigidity case, we pick a set of primes~$R$
  with $\abs{R} > 4t$ and each prime in $R$ larger than $4t$.  Then,
  for any subset $R'$ of $R$ of size at least~$t$,
  $\prod_{p \in R'} p > B_\loc$.
  Indeed,
  \begin{equation}
    \sum_{p\in R'} \ln p - \ln B_\loc
     \ge t\ln (4t) - t(\ln t + \ln 2)
     = t \ln 2 > 0.
  \end{equation}
  Therefore if we pick a prime~$p$ randomly from~$R$,
  $P_\loc$ is non-zero modulo~$p$ with probability at least $3/4$, and the
  Schwartz-Zippel Lemma tells us that if $P_\loc$ is non-zero
  modulo~$p$ we (correctly) answer ``yes''
  with probability at least $3/4$.  We conclude that our total
  probability of returning the correct answer is at least $\frac{3}{4}
  \cdot \frac{3}{4} > \frac{1}{2}$, as desired.

  Similarly, for global rigidity, we find
  \begin{equation}
    \begin{split}
      \ln B_\glob \le \ln B_1 + \ln B_2 
        &\le (s+1) e_\Gamma \ln e_\Gamma + (s+1) e_\Gamma \ln 2 + s\ln s\\
        &< v e_\Gamma \ln e_\Gamma + ve_\Gamma \ln 2 + v \ln v\\
        &< v e_\Gamma \ln (4 v e_\Gamma).
    \end{split}
  \end{equation}
  Again, suppose we pick a set of primes~$R$ with $\abs{R}> 4 v
  e_\Gamma$ and each prime bigger than $4 v e_\Gamma$.  Then for any
  $R' \subset R$ with $\abs{R'} > \abs{R}/4$, we have $\prod_{p \in
    R'} p < B_\glob$, as desired.
\end{proof}

%%% Local Variables: 
%%% mode: latex
%%% TeX-master: "main"
%%% End: 

\section{Smooth higher-dimensional flexes}
\label{sec:flexes}

Let us now turn to the issue of \emph{higher-dimensional flexing}.
Given incongruent frameworks~$\rho$ and~$\sigma$ in~$\EE^d$ with the
same edge lengths, how large must $a$
be so that $\rho$ and $\sigma$
are connected by a smooth path of frameworks in $\EE^{d+a}$ with
constant edge lengths?  (Such a path is called a \emph{flex} in
$\EE^{d+a}$.)
Bezdek and Connelly~\cite{BC04:KneserPoulsen}
have shown that an arbitrary pair of frameworks in~$\EE^d$  with the
same edge
lengths can be connected by a smooth path in~$\EE^{2d}$ (so
$a=d$ always suffices),
while Belk and Connelly~\cite{BC07:RigiditySimplexFlaps} exhibited a
$d$-dimensional framework of  the ``$d$-simplex
with flaps'' which is not globally rigid in~$\EE^d$ but still locally
rigid all the way up to $\EE^{2d-1}$ (so $a=d$
can be necessary).  Thus for
arbitrary frameworks one may have to go up to twice the
dimension to get any flexibility at all.  Theorem~\ref{thm:flex},
which we will now prove, states that
the situation is rather different for generic frameworks: 
for a generic framework~$\rho$ which is not globally rigid, there is
some other incongruent framework~$\sigma$ so that $\rho$ can be
flexed to $\sigma$ in
$\EE^{d+1}$ (so $a=1$ suffices).
The question of how large $a$ must be
to reach \emph{every} alternative framework with
the same edge lengths as a given generic one remains open.

The rest of this section is devoted to the proof of Theorem~\ref{thm:flex}.
Note that for the conditions of the theorem to be satisfied $\Gamma$
must have at
least $d+2$ vertices and must not have a minimal stress kernel.

\begin{definition}
Given a framework $\rho \in C^d(\Gamma)$,
and a stress $\Omega \in S(\rho)$,
  the space of \emph{lifted stress satisfiers} 
$\tA(\Omega)$ is the space of all
  ($d+1$)-dimensional frameworks of $\Gamma$
 that satisfy $\Omega$.
Also let $\trh$ be the framework in $\tA(\Omega)$
defined by embedding
$\rho$ in the first $d$ dimensions of $\EE^{d+1}$.
In our proof,
our $(d+1)$-dimensional flex will in fact stay in the
space $\tA(\Omega)$.
Because $\tA(\Omega)$ is isomorphic to
$K(\Omega)^{d+1}$, it is a ($kd+k$)-dimensional
linear space, where $k := \dim(K(\Omega))$.
  Finally, let $\tl$ be the length squared map on $\tA(\Omega)$
  and define $\tB(\Omega) \coloneqq \tl(\tA(\Omega))$.
\end{definition}

We
summarize here the main spaces and maps
that will be used in this proof.
\[
\begin{tikzpicture}
  \matrix(spaces) [matrix of math nodes, column sep=1cm, row sep=1cm]
  {
F(\rho,\Omega) & \tA(\Omega) & \tB(\Omega)  & L(\Omega) & \RR^e \\
        & F(\rho,\Omega)/\Eucl(d) & K(\Omega)& C^1(\Gamma)\\
};

\draw[right hook->] (spaces-1-1)--(spaces-1-2);
\draw[->>] (spaces-1-2) -- 
node[cdarrow]{\text{$\tl$}} 
(spaces-1-3);
\draw[right hook->] (spaces-1-3)--
node[cdarrow]{\text{same}} node[cdarrow,below]{\text{dim}} 
(spaces-1-4);
\draw[right hook->] (spaces-1-4)--(spaces-1-5);
\draw[->>] (spaces-1-1)--(spaces-2-2.north west);
\draw[->>] (spaces-2-2)--
node[cdarrow]{\text{$\varphi$}}
(spaces-2-3);
\draw[right hook->] (spaces-2-3)--(spaces-2-4);

\end{tikzpicture}
\]

Let $F(\rho,\Omega)$ be $\tl^{-1}(\ell(\rho))$, the fiber of $\tl$
over the point $\ell(\rho)$, and
let $F_0(\rho,\Omega)$ be the connected component containing~$\trh$.
Thus $F(\rho,\Omega)$ consists of 
points in $\tA(\Omega)$ with
the same edge lengths as $\rho$. 
We will first show that for any generic $\rho$, 
and suitably generic $\Omega \in S(\rho)$,
$F(\rho,\Omega)$ is a smooth
manifold. To do this, we will
need to know that $\ell(\rho)$ is a regular value 
of $\tl$ (Lemma~\ref{lem:flex-reg}). To compute the dimension of
$F(\rho,\Omega)$, we will
show that $\tB(\Omega)$ is the same dimension as 
the space $L(\Omega)$, which was defined in Definition~\ref{def:BL}.

By definition, $\trh$ 
has all zeros in its last ($(d+1)$'st) coordinate.
Thus our plan is to show that
there is another incongruent
framework in $F_0(\rho,\Omega)$
that also has all zeros 
in its last coordinate and is thus in $C^d(\Gamma)$.
To deal with congruences, we will mod out
by the Euclidean transformations in the first $d$ coordinates.
When $\Gamma$ does not have a minimal stress kernel, 
the
singularities of $F_0(\rho,\Omega)/\Eucl(d)$  
will be of codimension 2 or greater
(Lemma~\ref{lem:flex-codim}; compare Proposition~\ref{prop:sing-codim-Arho}).
We then consider a proper map $\varphi$ from
$F_0(\rho,\Omega)/\Eucl(d)$
that looks at the last coordinate of 
a framework, so
$[\trh]$ maps to~$0$.
We  will prove
that this map has a well defined and even degree.
This will 
prove the existence of a second
framework $\ts$ in $F_0(\rho,\Omega)$ which maps to zero under
$\varphi$
and thus comes from another framework~$\sigma$
in $C^d(\Gamma)$ with the same edge lengths.

Since $\trh$ and $\ts$ are both in $F_0(\rho,\Omega)$, which is connected by
assumption,
there must be a smooth path within $F(\rho,\Omega)$ which 
connects the two frameworks, which is the desired smooth flex
in~$\EE^{d+1}$.

\medskip

We now proceed to fill in the details.
In order to understand the fibers of $\tl$,
we would like to apply
Sard's Theorem, but we first need to know
that $\ell(\rho)$ is a generic point in~$L(\Omega)$ in an appropriate
sense. 
Here we use the standard notation $\QQ(x)$, where $x\in\RR^n$, to mean
the smallest field containing $\QQ$ and all the coordinates of~$x$.

\begin{lemma}\label{lem:generic-in-A}
  Let $\rho \in C^d(\Gamma)$ and $\omega \in S(\rho)$ such that 
  $(\ell(\rho),\omega)$ is generic in $C_\oM$.  Then 
 $\ell(\rho)$ is generic inside
  $L(\Omega)$, where we consider
  $L(\Omega)$ to be defined over $\QQ(\Omega)$.
\end{lemma}

\begin{proof}
Suppose that $\ell(\rho)$ is not generic in $B(\Omega)$ over
  $\QQ(\Omega)$.  Then, by the definition of genericity, there is some
  polynomial function~$P$ with coefficients in $\QQ$ so that
  $P(\ell(\rho),\Omega) = 0$ while there is some $\rho' \in A(\Omega)$
  so that $P(\ell(\rho'),\Omega) \ne 0$.  We can suppose that
  $d\ell_{\rho'}$ has maximal rank and 
  $\ell(\rho')$ is smooth in $\oM$ 
since 
having non maximal rank,
having an image under $\ell$ that is non smooth in $\oM$,
and the vanishing of $P$
are all defined by algebraic equations
  that do not vanish identically on $A(\Omega)$.  Since 
$\Omega$ is an equilibrium stress for $\rho'$,
$d\ell_{\rho'}$ has maximal rank,
and $\ell(\rho')$
is smooth,
from Lemma~\ref{lem:tangent-stress}
$(\ell(\rho'),
  \omega)$ must be in $C_\oM$. Thus $P$ is a polynomial function that
  does not vanish identically over $C_\oM$, contradicting the
  genericity of $(\ell(\rho), \omega)$ in $C_\oM$.

  Since $\ell(\rho)$ is generic in $B(\Omega)$, it is also generic 
in its Zariski closure $L(\Omega)$.
\end{proof}

Next we can look at the  map  $\tl$ and its fibers.  
Let $F(\rho,\Omega) \coloneqq
\tl^{-1}(\ell(\rho))$ be the fiber of $\tl$ at the point $\ell(\rho)$.
Also let $F_0(\rho,\Omega)$ be the component of 
$F(\rho,\Omega)$ 
that includes~$\rho$.

\begin{lemma}\label{lem:flex-reg}
Suppose
$(\ell(\rho),\omega)$ is generic in~$C_\oM$.
Then $F(\rho,\Omega)$ 
is a smooth manifold.
\end{lemma}

\begin{proof}
From Lemma~\ref{lem:generic-in-A}, $\ell(\rho)$ is a generic point
in $L(\Omega)$.
Thus from 
the algebraic version of Sard's
Theorem over the field $\QQ(\Omega)$, $\ell(\rho)$ is a regular
value of~$\tl$.
Therefore, from the implicit function theorem, 
$F(\rho,\Omega)$ 
is a smooth manifold.
\end{proof}

\begin{lemma}
  \label{lem:tB-dim}
Suppose that $\Gamma$, 
a graph with $d+2$ or more vertices,
is generically locally rigid in $\EE^d$ 
and $(\ell(\rho),\omega)$ is generic in~$C_\oM$. 
  Then the semi-algebraic set $\tB(\Omega)$ has the same dimension as
  $B(\Omega)$ and $L(\Omega)$, namely $kd - \binom{d+1}{2}$.
\end{lemma}

\begin{proof}
In one direction, $\tB(\Omega) \superset B(\Omega)$
and so clearly is at least of dimension
$kd-\binom{d+1}{2}$.
Recall from Proposition~\ref{prop:Bflat} that $B(\Omega)$ is a flat
space and so  is contained in
$L(\Omega)$, a linear space of the same dimension.
In the other direction, we need to show that 
$\tB(\Omega)$ is contained in 
$L(\Omega)$.
As described in Lemma~\ref{lem:B-chord},
$B(\Omega)$
coincides with 
$\chord^d(\ell(K(\Omega)))$
and
similarly $\tB(\Omega)$ is $\chord^{d+1}(\ell(K(\Omega)))$.
But then $\tB(\Omega) \subset \chord^{2d}(\ell(K(\Omega))) =
\chord^2(B(\Omega)) \subset L(\Omega)$, as desired.
\end{proof}

\begin{corollary}
  \label{cor:inL}
Suppose that $\Gamma$, 
a graph with $d+2$ or more vertices,
is generically locally rigid in $\EE^d$ 
and $(\ell(\rho),\omega)$ is generic in~$C_\oM$.
Then $F(\rho,\Omega)$ 
is a smooth manifold of dimension $k+\binom{d+1}{2}$.
\end{corollary}

\begin{proof}
$F(\rho,\Omega)$ is a smooth manifold by Lemma~\ref{lem:flex-reg}.  To find
its dimension,
subtract the dimension of the image
from the dimension of the domain:
\begin{equation}
\textstyle
\dim(F(\rho,\Omega)) = \dim \tA(\Omega) - \dim \tB(\Omega)
  = k(d+1) - \left(kd - \binom{d+1}{2}\right) = k+\binom{d+1}{2}.
\qedhere
\end{equation}
\end{proof}

Next we mod out by the the group consisting of Euclidean
transformations on the first $d$ coordinates. (We do not mod out by
Euclidean transforms involving the last coordinate, as we want to be
able to detect when the framework lies completely in $\EE^d$.)  This gives us
the quotient 
$F(\rho,\Omega)/\Eucl(d)$.

\begin{lemma}\label{lem:flex-codim}
Suppose that $\Gamma$, 
a graph with $d+2$ or more vertices,
is generically locally rigid in $\EE^d$,
$(\ell(\rho),\Omega)$ is generic in~$C_\oM$,
and $k>d+1$.
Then $F(\rho, \Omega)/\Eucl(d)$
is a smooth stratified space of dimension~$k$ with
singularities of codimension at least~$2$.
\end{lemma}

\begin{proof}
From Corollary~\ref{cor:inL}, $F(\rho,\Omega)$ is a smooth manifold.
So by the stratified structure on quotient spaces as in
Lemma~\ref{lem:sing-codim-gen}, 
$F(\rho,\Omega)/\Eucl(d)$
is a smooth stratified space with
singularities at quotients with non\hyp trivial stabilizer, i.e.,
frameworks
whose projection down to the first $d$ coordinates
span a proper affine subspace of $\EE^d$.  

First we argue that 
such singularities occur only at frameworks 
that span exactly
a $d$\hyp dimensional subspace of $\EE^{d+1}$ 
and that project down to
subspace of dimension exactly $d-1$ in the first $d$ coordinates.  
In particular, such singularities cannot have an affine span of
less than $d$ in $\EE^{d+1}$.
Any framework with a smaller affine span would give a framework in
$\EE^{d'}$ for some $d' < d$.  Considered as a framework in $\EE^d$,
this would mean that $\ell(\rho)$ is
not a regular value of $\ell$, contradicting
Lemma~\ref{lem:regval}.  Thus the only possibility for a
larger-than-expected stabilizer is a framework with $d$-dimensional
span whose projection is smaller than $d$-dimensional, as claimed.

The
dimension of these singular frameworks is  $d +
\binom{d+1}{2}$: it is $d$ (the dimension of the space of hyperplanes
in $\EE^d$, the possible choices of $(d-1)$-dimensional hyperplanes)
plus $\binom{d+1}{2}$ (the dimension of $\Eucl(d)$, the possible
choices of the framework within a single hyperplane; here
we use local rigidity in~$\EE^d$ and the fact that 
$\ell(\rho)$ is a regular value of $\ell$).  The singular frameworks are
therefore codimension $k-d$ inside $F(\rho,\Omega)$. 

Since the stabilizer of
these singular frameworks is always $O(1) \isom \ZZ/2$, which is
discrete, the codimension of the singular set does not change inside
the quotient.  In particular, since $\Gamma$ does not have a minimal
stress kernel, the codimension is always at least~$2$.

Since generically the stabilizer is trivial, the dimension of the
quotient $F(\rho,\Omega)/\Eucl(d)$ is
$\dim(F(\rho,\Omega)) - \dim(\Eucl(d)) = k$.
\end{proof}

Now we will look at the last coordinate of frameworks in 
$F_0(\rho,\Omega)/\Eucl(d)$.
Since all such frameworks are in $\tA(\Omega)$, this coordinate must
be in $K(\Omega)$. 
(Here we are thinking of $K(\Omega)$ as a space of
frameworks in $\EE^1$.)
Points that map to $0$ under this map represent
classes of frameworks 
that lie entirely in the space spanned by 
the first
$d$ coordinates in $\EE^{d+1}$.  Like all frameworks in $F(\rho,
\Omega)$, they have the same edge lengths
as $\rho$.

\begin{definition}
Let the map~$\varphi$ from 
$F_0(\rho,\Omega)/\Eucl(d)$
to $K(\Omega)$ be given
by simply 
looking at the last coordinate.
\end{definition}

\begin{lemma}
  \label{lem:phi-proper}
  The map $\varphi: F_0(\rho,\Omega)/\Eucl(d) \to K(\Omega)$ is proper.
\end{lemma}

\begin{proof}
  Think of $F_0(\rho,\Omega)/\Eucl(d)$ as a subset of the product
  $C^d(\Gamma)/\Eucl(d) \times C^1(\Gamma)$.  It suffices to show that
  $\varphi^{-1}(P)$ is bounded for any compact subset~$P$ of
  $K(\Omega)$.  (It is automatic that $\varphi^{-1}(P)$ is closed.)
  Inside $C^d(\Gamma)/\Eucl(d)$, all of $F_0(\rho,\Omega)/\Eucl(d)$ 
is bounded (as the
  edge lengths are bounded by the edge lengths of~$\rho$).  On the
  other hand, since~$P$ is a bounded subset of $K(\Omega)$, the
  projection of $\varphi^{-1}(P)$ to $C^1(\Gamma)$, which is 
just $P$ itself, is by supposition bounded.
\end{proof}
(Compare Lemma~\ref{lem:lisproper}.)

And now we are in position to complete the proof of our Theorem.

\begin{proof}[Proof of Theorem~\ref{thm:flex}]
If $\Gamma$ is not generically locally rigid in~$\EE^d$, then clearly
generic frameworks can be flexed in $\EE^{d+1}$ (as they can be flexed
in $\EE^d$).  So from now on we assume~$\Gamma$ is generically locally
rigid in $\EE^d$.

By Lemma~\ref{lem:generic-omega},
for any generic $\rho$, there is
an $\omega \in S(\rho)$ such that $(\ell(\rho),\omega)$ is generic in 
$C_\oM$.
Choose this $\omega$
to define $\varphi$. Since $\Gamma$ does not have a minimal stress kernel,
$k \geq k_{\min} > d+1$. (In fact for such a generic $\Omega$,
$k$ must equal $k_{\min}$).

By Lemma~\ref{lem:flex-codim}, when $k>d+1$ 
the space $F_0(\rho,\Omega)/\Eucl(d)$ has singularities
of high codimension.
By Lemma~\ref{lem:phi-proper}, $\varphi$ is proper, and by
Lemma~\ref{lem:flex-codim} again, the domain of $\varphi$ has the same
dimension
as its range (which is $k$).
Thus
$\varphi$ 
has  a well-defined mod-two degree by Corollary~\ref{cor:degree-sing}.
The frameworks in $F_0(\rho,\Omega)$ have fixed edge lengths, which implies
that the image frameworks in $K(\Omega)$, considered as a subset of
$C^1(\Gamma)$, have all points within a bounded distance of each
other.  In particular $\varphi$ is not onto,
and so its mod-two degree must be  zero.
The preimages of~$0$ are congruence classes of frameworks of $\Gamma$ 
that lie in $\EE^d$ and have the
same edge lengths as~$\rho$.

Next, by the analysis in Lemma~\ref{lem:flex-codim}, the point~$0$
cannot be the image of a singularity of $F_0(\rho,\Omega)/\Eucl(d)$.
Furthermore, $0$ is a regular
value of~$\varphi$: An element of the kernel of $d\varphi$
at $[\sigma]$ for some $\sigma$ in the inverse image of~$0$ is an
infinitesimal $d$-dimensional flex of the framework ~$\sigma$ in
$C^d(\Gamma)$.  But
Lemma~\ref{lem:regval} tells us that every
framework
in $C^d(\Gamma)$
with the same edge lengths as~$\rho$ is infinitesimally rigid.
Thus there are an even number of points in $\varphi^{-1}(0)$.  Let
$[\ts]$ be another such point.
Since $\trh$ and $\ts$ are points in a connected smooth manifold, there is a
smooth path connecting them.  This smooth path is the desired path of
frameworks of $\Gamma$ in
$\EE^{d+1}$ with
constant edge lengths.
\end{proof}

%%% Local Variables: 
%%% mode: latex
%%% TeX-master: "main"
%%% End: 

\bibliographystyle{hamsplain}
\bibliography{graphs,alggeom,topo}

\end{document}